# ON THE 2D ISING WULFF CRYSTAL NEAR CRITICALITY

By R. Cerf and R. J. Messikh

*Université Paris-Sud*

We study the behavior of the two-dimensional Ising model in a finite box at temperatures that are below, but very close to, the critical temperature. In a regime where the temperature approaches the critical point and, simultaneously, the size of the box grows fast enough, we establish a large deviation principle that proves the appearance of a round Wulff crystal.

**1. Introduction.** The Ising model in two dimensions is the first model where phase transition and non mean-field critical behavior has been established by Onsager [32] in 1944. It is also for that particular model that phase coexistence was rigorously studied and led to the first microscopic justification of the Wulff crystal. The first proof by Dobrushin, Kotecký and Shlosman [18] is valid for temperatures that are much lower than the critical point. Simplifications by Pfister [33], Ioffe and Schonmann [24–26] improved the result up to the critical point. These results in two dimensions rely on the study of contours to analyze large deviations of surface order. The extension of the Wulff construction to the Ising model in dimensions greater or equal to three required new techniques such as block coarse graining and the use of tools coming from geometric measure theory. This was achieved by Cerf and Pisztora [12] and Bodineau [6]. These results were initially valid up to the slab "percolation" threshold, and recently Bodineau [7] proved that this threshold is indeed the usual critical point thus extending the results of [12, 35] up to the critical temperature $T_c$. A two-dimensional analogue of the coarse graining developed in [35] is the subject of [15], thereby providing a unified approach to treat the problem for all dimensions. The appearance of the Wulff crystal has been proved in other "percolation"-type models as well, for example, [3, 4] in two dimensions and [9, 13] for dimensions greater or equal to three.

In all the works described above, the temperature has been kept fixed away from the critical temperature. Our main goal is to study the impact of the presence of a second-order phase transition on the phase coexistence phenomenon. We do this by analyzing phase coexistence in a regime where the temperature approaches the critical point from below while simultaneously taking the thermodynamical limit.

A priori, one can expect that a second-order phase transition has a non trivial effect on phase coexistence. Indeed, when approaching the critical point, the basic quantities describing the model either diverge or stay finite but have divergent derivatives. In the second case, they decay as a power law giving rise to critical exponents. This critical behavior and the existence of these exponents is conjectured for a wide family of two-dimensional statistical mechanics models. The existence of some of these critical exponents in a strong sense is an important ingredient in our analysis. This ingredient is available for the two-dimensional Ising model. For this particular model, the relevant statistical mechanics quantities can be computed explicitly, giving rise to beautiful identities such as Onsager's formula for the magnetization [32] and an explicit formula for the surface tension that describes the geometry of the Wulff droplet in terms of random walks [30, 31]. Such results can be obtained, for example, by using the dimer representation of Kasteleyn [27].

The probabilistic understanding of the critical phenomena is a very active field nowadays. In the case of Bernoulli site percolation on the planar triangular lattice, the existence and the identification of the critical exponents have been rigorously established by Smirnov and Werner [38] after the groundbreaking work of Schramm [36] and Smirnov [37]. In [38], the existence of the critical exponents has been explained in a probabilistic manner. Indeed, this work establishes a rigorous link with the conformal invariance of the scaling limit of critical percolation described by the Schramm–Loewner evolution process. Regarding the above mentioned results, the reader may wonder why we do not investigate the influence of a phase transition on phase coexistence in the a priori simpler and better understood *independent* Bernoulli percolation model instead of the *dependent* spins of the Ising model. The reason is that despite the spectacular progress in the understanding of criticality of independent percolation, essential properties of the critical exponents are still unaccessible by other methods than explicit computation (see in particular the open question 3 at the end of [38]). And since explicit computations work only for the two-dimensional Ising model, our results are confined to this particular model.

**2. Statement of the main results.** Consider the Ising model at temperature $T < T_c$, defined on a square box $\Lambda(n)$ of side length $n \in \mathbb{N} \setminus \{0\}$ and submitted to plus boundary conditions. For every spin configuration $\sigma \in \{-1, 1\}^{\Lambda(n)}$, we associate a random signed measure $\sigma_{n,T}$ on the unit square $Q = [-1/2, 1/2]^2 \subset \mathbb{R}^2$ defined by

$$\sigma_{n,T} = \frac{1}{m^*(T)n^2} \sum_{x \in \Lambda(n)} \sigma(x)\delta_{x/n},$$

where $\delta_{x/n}$ is the Dirac mass at $x/n$, and $m^*(T)$ is the spontaneous magnetization at temperature $T$. The expectation $b_n$ under $\sigma_{n,T}$ is

$$b_n = \int_Q x \, d\sigma_{n,T}(x) = \frac{1}{m^*(T)n^3} \sum_{x \in \Lambda(n)} \sigma(x)x.$$

The main result of our paper is the following convergence theorem for $\sigma_{n,T}$ under a conditioned Ising measure:

THEOREM 1. *Let $0 < \delta < \pi$. Let $B(\delta)$ be the ball of radius $\sqrt{\delta/\pi}$ and $w_n$ be the random measure defined by*

$$w_n(x) \, dx = \left(1_Q(x) - 2 \cdot 1_{B(\delta)}\left(\frac{b_n}{2} + x\right)\right) dx.$$

*This is the measure having density $-1$ on $B(\delta) - b_n/2$ and $1$ on the complement. Under the conditional probability,*

$$\mu_{\Lambda(n)}(\cdot) = \mu_{\Lambda(n)}^{+,T}\left(\cdot \left| \frac{1}{n^2} \sum_{x \in \Lambda(n)} \sigma(x) \leq (1-\delta)m^*(T) \right.\right),$$

*the difference between the random measures $\sigma_n$ and $w_n$ converges weakly in probability towards $0$ when $n \uparrow \infty$ and $T \uparrow T_c$ in such a way that $n(T_c - T)^{20} \uparrow \infty$, and $\log(n)/\log(1/(T_c - T))$ stays bounded. That is, for any continuous function $f: Q \to \mathbb{R}$,*

$$(1) \qquad \forall \varepsilon > 0 \qquad \lim_{n,T} \mu_{\Lambda(n)}(|\sigma_n(f) - w_n(f)| \geq \varepsilon) = 0.$$

*The probabilities of the deviations are of order $\exp(-\text{constant}(T_c - T)n)$.*

The last sentence of the theorem means the following. For any continuous function $f: Q \to \mathbb{R}$, any $\varepsilon > 0$, there exist positive constants $b, c$ depending on $f, \varepsilon$ such that

$$\mu_{\Lambda(n)}\left(\left| \frac{1}{m^*(T)n^2} \sum_{x \in \Lambda(n)} \sigma(x)f\left(\frac{x}{n}\right) \right.\right.$$

$$\left.\left. + 2\int_{B(\delta)} f\left(-\frac{b_n}{\delta} + x\right) dx - \int_Q f(x) \, dx \right| > \delta\right)$$

$$\leq b \exp(-c(T_c - T)n).$$

The main assertion of the theorem is that conditioned on having a defect of magnetization, the random measure $\sigma_{n,T}$ looks like a measure whose density is an indicator of the Wulff crystal which turns out to be an ordinary circle near the critical point. In other words, the defect of magnetization concentrates into a circular region. Note that in our regime, the shape of the Wulff crystal is no more affected by the geometry of the square lattice.

Theorem 1 is a consequence of De Giorgi's isoperimetric inequality [16] and a large deviation principle (LDP) that we prove in this paper. The assumption on $\log(n)/\log(1/(T_c - T))$ is a side hypothesis. Although we believe that a proof in the case where this quantity diverges is possible using the ideas of the current paper, we could not find a proof that includes both cases. Since we would like to study regimes that are as close as possible to criticality, we decided to treat the case where $\log(n)/\log(1/(T_c - T))$ stays bounded. The exponent 20 in the statement of the theorem is not optimal. Indeed, if we introduce the quantity

$$\nu_{\mathcal{W}} = \inf \left\{ \begin{array}{l} \gamma > 0 \text{ such that the convergence (1) is valid when } n \uparrow \infty \\ \text{and } T \uparrow T_c \text{ in such a way that } n(T - T_c)^\gamma \uparrow \infty \end{array} \right\},$$

then our result states that $\nu_{\mathcal{W}} \leq 20$. We believe that $\nu_{\mathcal{W}} = 1$, that is, it should be the critical exponent for the correlation length. For percolation-type models we can introduce a similar exponent that characterizes the maximal regime where a Wulff droplet near criticality appears. We believe that in the case of the two-dimensional Bernoulli percolation $\nu_{\mathcal{W}} = 4/3$.

For the two-dimensional Ising model, several difficulties have to be overcome to obtain the right exponent. In a heuristic manner, what prevents us to go from $\nu_{\mathcal{W}} \leq 20$ to $\nu_{\mathcal{W}} \leq 5$ is the lack of the van den Berg–Kesten inequality for the dependent random cluster model. Then to go from $\nu_{\mathcal{W}} \leq 5$ to $\nu_{\mathcal{W}} \leq 2$, one has to have a better understanding of the influence of the boundary conditions when approaching criticality. More precisely, one has to understand better how weak mixing properties in the sense of [1] behave close to the critical point. The gap $\nu_{\mathcal{W}} \leq 2$ to $\nu_{\mathcal{W}} \leq 1$ is related to the speed of convergence of the empirical magnetization to its thermodynamical limit in a regime where we approach the critical point.

An alternative approach to investigate the identity $\nu_{\mathcal{W}} = 1$ would be to use the finer DKS-theory along the work of [18, 26, 33, 34]. Indeed, it appears from the comments of one of the referees that such an approach might work. Here we chose to use an approach along the lines of [9, 35]. We believe that such an approach is interesting by itself. Its advantage is its potential application in studying models in higher dimensions. The study of phase coexistence near criticality in higher dimensions, say in the simpler case of Bernoulli percolation, would require some information on the critical behavior of the surface tension. To our knowledge such results are, at the present time, not available. If we suppose that a surface tension near criticality exists in higher dimensions, all our arguments are likely to work in this case except for the interface lemma. Indeed, our analysis shows that this part of the proof would need a nontrivial adaptation. We will give more details about this problem in Section 5.4.

The computation of an explicit bound for $\nu_{\mathcal{W}}$ would not have been realized without the kind advice of Charles Pfister. We thank him warmly for drawing our attention to the paper [8].

2.1. *The large deviation principle.* Next, we give the notions necessary to state our LDP that describes phase coexistence near the critical point.

2.1.1. *Sets of finite perimeter.* A fundamental quantity in our study is the perimeter of a set. In order to prove our LDP, we need to define this quantity for Borel subsets of $Q = [-1/2, 1/2]^2$ that may have irregular boundaries. We give next the distributional definition of the perimeter.

DEFINITION 2. Let $A$ be a Borel subset of $Q$, its perimeter is defined as

$$\mathcal{P}(A) = \sup \left\{ \int_A \operatorname{div} f(x) \, d\mathcal{L}^2(x) : f \in C_c^\infty(\mathbb{R}^2, B(0,1)) \right\},$$

where $\mathcal{L}^2$ is the Lebesgue measure on $\mathbb{R}^2$, $C_c^\infty(\mathbb{R}^2, B(0,1))$ is the set of $C^\infty$ vector fields from $\mathbb{R}^2$ to the Euclidean unit ball $B(0,1)$ having a compact support and div is the usual divergence operator. The set $A$ is said to have finite perimeter if $\mathcal{P}(A)$ is finite.

If the boundary of $A$ is smooth then an application of the Gauss–Green theorem gives that $\mathcal{P}(A) = \mathcal{H}(\partial A)$ where $\mathcal{H}$ is the one-dimensional Hausdorff measure. We denote by $\mathcal{M}(Q)$ the vector space of the finite signed Borel measures on $Q$. We equip $\mathcal{M}(Q)$ with the weak topology, that is, the coarsest topology for which the linear functionals

$$\nu \in \mathcal{M}(Q) \mapsto \int f \, d\nu, \qquad f \in C_c(\mathbb{R}^2, \mathbb{R}),$$

are continuous where $C_c(\mathbb{R}^2, \mathbb{R})$ is the set of the continuous functions from $\mathbb{R}^2$ to $\mathbb{R}$ having compact support. The rate function of our LDP is the map

$$\mathcal{J} : \mathcal{M}(Q) \to [0, +\infty],$$

$$\nu \mapsto \mathcal{J}(\nu) = \begin{cases} \tau_c \mathcal{P}(A), & \text{if there exists a Borel set } A \text{ such that} \\ & \dfrac{d\nu}{d\mathcal{L}^2} = -1_A + 1_{Q \setminus A}, \\ \infty, & \text{otherwise.} \end{cases}$$

The positive constant $\tau_c$ will be defined later; it plays the role of the unit length by which we measure the perimeter.

THEOREM 3.    *If $T \uparrow T_c$ and $n \uparrow \infty$ in such a way that $n(T_c - T)^{20} \uparrow \infty$ and such that $\log(n)/\log(1/(T_c - T))$ stays bounded, then the law of the random measure $\sigma_{n,T}$ under $\mu_{\Lambda(n)}^{+,T}$ satisfies a large deviation principle in $\mathcal{M}(Q)$ with respect to the weak topology. The speed of the LDP is $(T_c - T)n$, and the good rate function is $\mathcal{J}$; that is, for any Borel subset $\mathbb{M}$ of $\mathcal{M}(Q)$,*

$$-\inf\{\mathcal{J}(\nu): \nu \in \overset{\circ}{\mathbb{M}}\} \leq \liminf_{n,T} \frac{1}{(T_c - T)n} \log \mu_{\Lambda(n)}^{+,T}[\sigma_{n,T} \in \mathbb{M}]$$

$$\leq \limsup_{n,T} \frac{1}{(T_c - T)n} \log \mu_{\Lambda(n)}^{+,T}[\sigma_{n,T} \in \mathbb{M}]$$

$$\leq -\inf\{\mathcal{J}(\nu): \nu \in \overline{\mathbb{M}}\}.$$

2.2. *Structure of the paper.* In Section 3, we define the Ising model and introduce the basic notions that we will use in the rest of the paper. In this section we also give the Ising-specific properties on which we rely. This way we isolate and minimize the use of the specificities of the Ising model. The other techniques are more robust and take their roots in the study of the Wulff crystal in dimensions greater or equal to three. Section 4 contains a decoupling lemma and preliminary block estimates. These results are based on the paper [11]. Section 5 contains the proof of the upper bound for our LDP. Finally Section 6 finishes the proof of the LDP by establishing the corresponding lower bound.

**3. Preliminaries.** In this section we define the Ising model and its representation in terms of the random cluster model. In a second part, we isolate the Ising-specific properties that are required for our study.

3.1. *The Ising model.* Let $Q = [-1/2, 1/2]^2$ be the centered planar unit box. For a positive integer $n$, we define the discrete set of sites $\Lambda(n) = nQ \cap \mathbb{Z}^2$ that we turn into a graph by considering the following set of edges:

$$\mathbb{E}(\Lambda(n)) = \{\{x, y\} \subset \Lambda(n): |x - y| = 1\},$$

where $|\cdot|$ is the usual Euclidean norm. We also define the boundary $\partial \Lambda(n)$ of the graph $\Lambda(n)$ by

$$\partial \Lambda(n) = \{x \in \Lambda(n): \exists y \in \Lambda(n)^c: |x - y| = 1\}.$$

For every value $\beta = 1/T > 0$ of the inverse temperature, the Ising model in $\Lambda(n)$ with $+1$ boundary conditions is the probability measure on the spin configurations $\Omega_{\Lambda(n)} = \{-1, +1\}^{\Lambda(n)}$ defined by

$$\forall \sigma \in \{+1, -1\}^{\Lambda(n)} \qquad \mu_{\Lambda(n)}^{+,\beta}[\sigma] = \frac{1}{Z_n^{+,\beta}} \exp(-\beta \mathcal{H}_n^+(\sigma)),$$

where

$$\mathcal{H}_n^+(\sigma) = -\frac{1}{2} \sum_{\substack{x,y \in \Lambda(n) \setminus \partial\Lambda(n) \\ |x-y|=1}} \sigma(x)\sigma(y) - \sum_{x \in \Lambda(n) \setminus \partial\Lambda(n)} \sum_{\substack{y \in \partial\Lambda(n) \\ |x-y|=1}} \sigma(x)$$

and $Z_n^{+,\beta}$ is the adequate normalization constant.

3.2. *The FK-representation.* There exists a useful and well-known coupling between the Ising model at inverse temperature $\beta$ and the random cluster model with parameter $q = 2$ and $p = 1 - \exp(-2\beta)$; see [19, 20]. We will use this coupling in order to derive several probabilistic estimates from the corresponding FK-percolation model. The coupling is a probability measure $\mathbb{P}_n^+$ on the edge–spin configuration space $\{0,1\}^{\mathbb{E}(\Lambda(n))} \times \{-1,+1\}^{\Lambda(n)}$. To construct $\mathbb{P}_n^+$ we first consider Bernoulli percolation of parameter $p$ on the edge space $\{0,1\}^{\mathbb{E}(\Lambda(n))}$; then we choose the spins of the sites in $\Lambda(n)$ independently with the uniform distribution on $\{-1,+1\}$, and finally we condition the edge–spin configuration on the event that there is no open edge in $\Lambda(n)$ between two sites with different spin values. The construction can be summed up with a formula; we have

$$\forall (\sigma, \omega) \in \{0,1\}^{\mathbb{E}(\Lambda(n))} \times \{-1,+1\}^{\Lambda(n)}$$

$$\mathbb{P}_n^+(\sigma, \omega) = \frac{1}{Z} \prod_{e \in \mathbb{E}(\Lambda(n))} p^{\omega(e)}(1-p)^{1-\omega(e)} \mathbb{1}_{(\sigma(x)-\sigma(y))\omega(e)=0},$$

where $Z$ is the adequate normalization constant. It can be verified that the marginal of $\mathbb{P}_n^+$ on the spin configurations is the Ising model with parameter $\beta$ given by the formula $p = 1 - \exp(-2\beta)$, and the marginal on the edge configurations is the random cluster measure with parameters $p$, $q = 2$ and subject to wired boundary conditions, that is, the probability measure on $\Omega_{\Lambda(n)} = \{0,1\}^{\mathbb{E}(\Lambda(n))}$ defined by

$$(2) \qquad \forall \omega \in \Omega_{\Lambda(n)} \qquad \Phi_{\Lambda(n)}^{p,w}[\omega] = \frac{1}{Z} q^{\text{cl}^w(\omega)} \prod_{e \in \mathbb{E}(\Lambda(n))} p^{\omega(e)}(1-p)^{1-\omega(e)},$$

where $\text{cl}^w(\omega)$ is the number of connected components with the convention that two clusters that touch the boundary $\partial\Lambda(n)$ are identified. This coupling says that one may obtain an Ising configuration by first drawing a FK-percolation configuration with the measure $\Phi_{\Lambda(n)}^{w,p}$, then coloring all the sites in the clusters that touch the boundary $\partial\Lambda(n)$ in $+1$ and finally coloring the remaining clusters independently in $+1$ and $-1$ with probability $1/2$ each. Also, the coupling permits one to obtain a $\Phi_{\Lambda(n)}^{w,p}$ percolation configuration by first drawing a spin configuration with $\mu_{\Lambda(n)}^{+,\beta}$, and then by declaring that all the edges between two sites with different spins are closed while the other edges are independently declared open with probability $p$ and closed with probability $1 - p$.

Let $\Lambda \subset \mathbb{Z}^2$ and $0 \le p \le 1$. In addition to the wired boundary conditions we will also work with *partially wired boundary conditions*. In order to define them, we consider a partition $\pi$ of

$$\partial \Lambda = \{x \in \Lambda : \exists y \in \mathbb{Z}^2 \setminus \Lambda, |x - y| = 1\}.$$

Let us say that $\pi$ consists of $\{B_1, \dots, B_k\}$ where the $B_i$ are nonempty disjoint subsets of $\partial \Lambda$ and such that $\bigcup_i B_i = \partial \Lambda$. For every configuration $\omega \in \Omega_\Lambda$, we define $\mathrm{cl}^\pi(\omega)$ as the number of open connected clusters in $\Lambda$ computed by identifying two clusters that are connected to the same set $B_i$. The $\pi$-wired FK-measure $\Phi_\Lambda^{p,\pi}$ is defined by substituting $\mathrm{cl}^\pi(\omega)$ for $\mathrm{cl}^w(\omega)$ in (2). We will denote the set of all partially-wired FK-measures in $\Lambda$ by $\mathcal{FK}(p, \Lambda)$. Note that $\Phi_\Lambda^{p,w}$ corresponds to $\pi = \{\partial \Lambda\}$. We define the FK-measure with free boundary conditions $\Phi_\Lambda^{p,f}$ as the partially-wired measure corresponding to $\pi = \varnothing$.

Let $U \subseteq V \subseteq \mathbb{Z}^2$. For every configuration $\omega \in \{0, 1\}^{\mathbb{E}(\mathbb{Z}^2)}$, we denote by $\omega_V$ the restriction of $\omega$ to $\Omega_V = \{0, 1\}^{\mathbb{E}(V)}$. More generally we will denote by $\omega_V^U$ the restriction of $\omega$ to $\Omega_V^U = \{0, 1\}^{\mathbb{E}(V) \setminus \mathbb{E}(U)}$. If $V = \mathbb{Z}^2$ or $U = \varnothing$, then we drop them from the notation. We will denote by $\mathcal{F}_V^U$ the $\sigma$-algebra generated by the finite-dimensional cylinders of $\Omega_V^U$. Note that every configuration $\eta \in \Omega_V$ induces a partially wired boundary condition $\pi(\eta)$ on the set $U$. The partition $\pi(\eta)$ is obtained by identifying the sites of $\partial U$ that are connected through an open path of $\eta^U$. We will denote by $\Phi_U^{p,\pi(\eta)}$ the corresponding FK-measure.

3.3. *Planar duality.* The duality of the FK-measures in dimension two is well known. In this paper we will use the notation of [15] that we summarize next. Let $0 \le p \le 1$ and $\Lambda$ be a connected subset of $\mathbb{Z}^2$. To construct the dual process of $\Phi_\Lambda^{p,w}$, we associate to each edge $e \in \mathbb{E}(\Lambda)$ the unique edge $\widehat{e}$ linking sites of $\mathbb{Z}^2 + (1/2, 1/2)$ which crosses orthogonally the edge $e$. The subgraph of $\mathbb{Z}^2 + (1/2, 1/2)$ consisting of the edges $\widehat{e}$ and their sites will be denoted by $(\widehat{\Lambda}, \mathbb{E}(\widehat{\Lambda}))$. To each configuration $\omega \in \Omega_\Lambda$ we associate the dual configuration $\widehat{\omega} \in \Omega_{\widehat{\Lambda}}$ defined by

$$\forall e \in \mathbb{E}(\Lambda) \qquad \widehat{\omega}(\widehat{e}) = 1 - \omega(e).$$

Similarly, for every event $A \subset \Omega_\Lambda$ we define the *dual event*

$$\widehat{A} = \{\eta \in \Omega_{\widehat{\Lambda}} : \exists \omega \in A, \widehat{\omega} = \eta\}.$$

The duality property asserts that

$$\Phi_\Lambda^{p,w}[A] = \Phi_{\widehat{\Lambda}}^{\widehat{p}, f}[\widehat{A}],$$

where $\widehat{p} = 2(1 - p)/(2 - p)$.

3.4. *The critical point.* It is known that the critical point of the Ising model on $\mathbb{Z}^2$ is given by the fixed point of a duality relation (see [22]). For the random cluster model with $q = 2$, the dual point $\widehat{p}$ is related to $p$ through the relation

(3)     $$\frac{p}{1-p}\frac{\widehat{p}}{1-\widehat{p}} = 2, \text{ and the fixed point is } p_c = \frac{\sqrt{2}}{1+\sqrt{2}}.$$

The duality relation (3) is translated in the Ising setting into

$$\sinh(2\beta)\sinh(2\widehat{\beta}) = 1, \text{ and the critical point is } \beta_c = \frac{\operatorname{arcsinh}(1)}{2}.$$

As we will work in the limit $p \to p_c$, it is worth noting that the derivative of the function $p \mapsto \widehat{p}$ is nonzero at $p_c$. Thus $p - p_c$ is of the same order as $p_c - \widehat{p}$ when $p \to p_c$. Also $\beta \mapsto p = 1 - \exp(-2\beta)$ has a nonzero derivative at $\beta_c$, and thus $p - p_c$ is of the same order as $\beta - \beta_c$ when $\beta \to \beta_c$.

For the general $q$-Potts model, the identification of the critical point and the self-dual point, that is, $p_c = \sqrt{q}/(1 + \sqrt{q})$, is still an open problem for the values $2 < q < 25$. When $q > 25.72$, this identity has been established, and in this situation the Potts model exhibits a first-order phase transition [21, 28]. Thus the 2d Ising model is the only two-dimensional Potts model exhibiting a second-order phase transition for which the critical point has been rigorously identified to be the self-dual point.

We end this section by setting the following convention concerning the use of the word *dual* in the rest of the paper: we always consider that the original model is the super-critical one, that is, $p > p_c$ which is defined on the edges of $\mathbb{Z}^2$. The *dual* model is always the dual of the super-critical model. That is, it is a sub-critical model defined on the edges of $\mathbb{Z}^2 + (1/2, 1/2)$ and at percolation parameter $\widehat{p} = 2(1-p)/(2-p) < p_c$. A dual path, circuit or site will always denote a path, circuit or site in $\mathbb{Z}^2 + (1/2, 1/2)$. The term *open dual* will always designate edges $\widehat{e}$ of $\mathbb{Z}^2 + (1/2, 1/2)$ that are open with respect to the dual configuration; that is, $\widehat{\omega}(\widehat{e}) = 1$. The law of the dual edges $\widehat{e}$ will always be the dual measure $\Phi^{\widehat{p}}$ which is sub-critical; that is, $\widehat{p} < p_c$.

3.5. *The surface tension and critical exponents.*

3.5.1. *The surface tension.* By duality, the surface tension $\tau_\beta$ of the two-dimensional Ising model at inverse temperature $\beta > \beta_c$ is given by the directional dependence of the exponential decay of the correlations at dual inverse temperature $\widehat{\beta} < \beta_c$;

$$\forall x \in \mathbb{Z}^2 \qquad \tau_\beta(x) = -\lim_{n\to\infty}\frac{1}{n}\log\mu_\infty^{\widehat{\beta}}[\sigma(0)\sigma(nx)] = -\lim_{n\to\infty}\frac{1}{n}\log\Phi_\infty^{\widehat{p}}[0 \leftrightarrow nx]$$

with $\widehat{p} = 1 - \exp(-2\widehat{\beta})$ and where we have used the FK-representation to derive the second equality. We will also consider the unique continuous extension of $\tau_\beta$ into a norm on $\mathbb{R}^2$.

In this paper, we are interested in the situation where the spatial scale $n$ goes to infinity, and simultaneously $\beta$ goes to $\beta_c$. To study phase coexistence in such a context, we first need to define the joint limit surface tension.

PROPOSITION 4.   *If $n \uparrow \infty$ and $\widehat{\beta} \uparrow \beta_c$ in such a way that $n(\beta_c - \widehat{\beta})/\log n \to \infty$, then uniformly on $x \in \mathbb{Z}^2$, we have that*

$$(4) \qquad \lim_{(n,\widehat{\beta}) \to (\infty, \beta_c)} -\frac{1}{(\beta_c - \widehat{\beta})n} \log \mu_\infty^{\widehat{\beta}}[\sigma(0)\sigma(nx)] = \tau_c|x|,$$

*where $\tau_c$ is a positive constant, and $|\cdot|$ is the Euclidean norm on $\mathbb{R}^2$.*

Note that $\tau_c$ does not depend on $x$ and that the appearance of the Euclidean norm means that the surface tension is isotropic near the critical point. The presence of the factor $(\beta - \beta_c)$ in the denominator of the limit is directly related to the critical exponent for the correlation length. In the planar Ising model this exponent is equal to 1.

The proof of the last proposition is an extension of the computation [14] of the *classical* surface tension at fixed temperature. Using subadditivity, it follows that

$$\forall x \in \mathbb{Z}^2 \qquad \log \mu_\infty^{\widehat{\beta}}[\sigma(0)\sigma(x)] \leq \tau_\beta(x).$$

Thus the upper bound of (4) follows directly from the formula obtained in [14] describing $\tau_\beta$. On the other hand, the corresponding lower bound is harder to obtain. It is obtained by extending the asymptotic computations of $\mu_\infty(\sigma(0)\sigma(nx))$ from a regime where $\widehat{\beta} < \beta_c$ is kept fixed and $n \uparrow \infty$ to a situation where $\widehat{\beta} \uparrow \beta_c$, and, simultaneously, $n \uparrow \infty$. This can be done by using, as in [14], Kasteleyn's dimer representation of the 2d Ising model [27]. These computations are long and rely on quite different mathematical tools, and thus we consider here (4) as a starting point and we will present the derivation of Proposition 4 in a separate paper [31]. In words, the extension [31] relies on a random walk interpretation of the computation [14]. This interpretation allows the derivation of Proposition 4 from classical moderate deviation results for the random walk. In fact, the isotropy of the right-hand side of (4) follows from the isotropy of the rate function for the moderate deviation principle of the simple random walk on $\mathbb{Z}^2$.

The regime in Proposition 4 is nearly optimal since, from results of [5], it appears that if $(n, \widehat{\beta}) \uparrow (\infty, \beta_c)$ in a regime where $n(\beta_c - \widehat{\beta})$ stays constant, then

$$\mu_\infty^{\widehat{\beta}}[\sigma(0)\sigma(ne_1)] \sim c\frac{1}{n^{1/4}},$$

where $e_1 = (1, 0) \in \mathbb{Z}^2$ and $c$ is a positive constant. Thus below the regime of Proposition 4 we do not expect phase coexistence to happen. In this situation the critical phenomena take over and hence we expect that $\nu_{\mathcal{W}} \geq 1$.

3.5.2. *The magnetization.* We will need to know how fast the magnetization is going to zero near the critical point. In our case, from Onsager's famous formula we know that

$$(5) \qquad \mu_{\infty}^{+,\beta}[\sigma(0)] \sim (\beta - \beta_c)^{1/8} \qquad \text{when } \beta \downarrow \beta_c.$$

The Ising model is the only Potts model for which this exponent has been computed. For independent site percolation on the planar triangular lattice, this exponent has been derived in [38] along with several other exponents. For the values $q = 1, 3, 4$, the existence of the exponent as well as its values are conjectured by physicists but currently not proved.

In addition to that, we shall need to estimate the speed of convergence of the empirical magnetization near the critical point. In order to do so, we rely on correlation inequalities that are specific to the Ising model. Furthermore, we rely on explicit computations to control the behavior of the relevant quantities near the critical point.

**4. Block arguments.** Besides the Ising-specific properties that we stated in the last section, our analysis is based on rather robust techniques that have been developed by Cerf and Pisztora [9, 10, 12, 13] to study phase coexistence in dimensions greater or equal to three. The probabilistic estimates are obtained by translating the relevant Ising events into the random cluster model via the FK-representation. In this paper, an essential tool in analyzing the random cluster model is an adaptation of block coarse graining techniques [35] to the situation where $p \downarrow p_c$.

4.1. *Notation and preparatory lemmas.* In this section, we introduce the notation used in coarse graining arguments and state useful preliminary estimates that we will use repeatedly in the rest of the paper.

4.1.1. *The rescaled lattice.* First we fix a positive integer $K$ that will typically depend on $p$ later on and will diverge when $p \downarrow p_c$. For each $\underline{x} \in \mathbb{Z}^2$, we define the block indexed by $\underline{x}$ as $B(\underline{x}) = \Lambda(K) + K\underline{x}$. Let $A$ be a region in $\mathbb{R}^2$. We define the rescaled region $\underline{A}$;

$$\underline{A} = \{\underline{x} \in \mathbb{Z}^2 : B(\underline{x}) \cap A \neq \varnothing\}.$$

From now on, underlining means that we are dealing with rescaled objects. For instance, $\underline{\Lambda}(n)$ means the rescaled box $\Lambda(n)$. Note that $|\underline{\Lambda}(n)|$ is now of order $n^2/K^2$ which is the order of the number of boxes necessary to cover $\Lambda(n)$.

4.1.2. *The lattice* $\mathbb{L}^\infty$. When dealing with block arguments it will be convenient to equip $\mathbb{Z}^2$ with another graph structure. We denote by $d_\infty$ the distance associated with the norm $|\cdot|_\infty$ defined by

$$\forall (x_1, x_2) \in \mathbb{R}^2 \qquad |(x_1, x_2)|_\infty = \max(|x_1|, |x_2|).$$

Thus $d_\infty(x, y) = |x - y|_\infty$. For every set $E \subseteq \mathbb{R}^2$ and positive real number $r$, we define the $r$-neighborhood of $E$ as

$$\mathcal{V}_\infty(E, r) = \{x \in \mathbb{R}^2 : d_\infty(x, E) < r\}.$$

We will also use the associated diameter given by

$$\mathrm{diam}_\infty(E) = \sup\{|x - y|_\infty : x, y \in E\}.$$

The new graph structure on $\mathbb{Z}^2$ is induced by the set of edges,

$$\mathbb{E}^{2,\infty} = \{\{x, y\} \subset \mathbb{Z}^2 : d_\infty(x, y) = 1\}.$$

Then the lattice $\mathbb{L}^\infty$ is the graph $(\mathbb{Z}^2, \mathbb{E}^{2,\infty})$. This lattice has the useful property that the exterior boundary of any connected finite set $A$ in $\mathbb{Z}^2$ is itself connected when regarded as a subgraph of $\mathbb{L}^\infty$; for a proof we refer the reader to [17].

4.1.3. *Block events.* For the renormalization to be useful it is almost always required to use block events on a set of blocks that are overlapping. Thus in addition to the partitioning blocks $B(\underline{x})$ we define the event blocks $B'(\underline{x})$ by setting

$$B'(\underline{x}) = \bigcup_{\substack{\underline{y} \in \mathbb{Z}^2 \\ |\underline{y} - \underline{x}|_\infty \leq 1}} B(\underline{y}).$$

4.1.4. *Rough estimates on the block process.* Given the events that describe a good block, we define the block process $(X(\underline{x}), \underline{x} \in \underline{\Lambda}(n))$ as the dependent site percolation process on $\underline{\Lambda}(n)$ that indicates if a block is good or not. We cite several rough estimates on the block process from [10]. The block process can be viewed as a dependent site percolation process where a site $\underline{x}$ is occupied if and only if $X(\underline{x}) = 0$. The occupied $\mathbb{L}^\infty$ cluster of the site $\underline{x}$, that is, the $\mathbb{L}^\infty$ connected component of the occupied sites containing $\underline{x}$, is then denoted by $\underline{C}(\underline{x})$. The next lemma is a standard counting Peierls argument:

LEMMA 5. *Suppose that there exists $\delta > 0$ such that*

$$\forall \underline{x} \in \mathbb{Z}^2 \qquad P[X(\underline{x}) = 0 | X(\underline{z}), |\underline{x} - \underline{z}|_\infty \geq 3] \leq \delta.$$

*There exists a constant $b$ such that, for any bounded open subset $O$ of $\mathbb{R}^2$, any $s, t > 0$, any $K, n \in \mathbb{N}$ with $n \geq K$,*

$$P[|\{\underline{x} \in \mathbb{Z}^2 : B(\underline{x}) \cap O \neq \varnothing, |\underline{C}(\underline{x})| \geq t\}| \geq s]$$

$$\leq 2 \sum_{j \geq s} \exp j\left(\frac{1}{t} \log \mathcal{L}^2(\mathcal{V}(O, 2)) + \log b + \frac{1}{9} \log \delta\right),$$

*where $\mathcal{V}(O, 2) = \{x \in \mathbb{R}^2 : d(x, O) \leq 2\}$.*

Here is the last rough estimate:

Lemma 6. *We consider the box $\Lambda(n)$ rescaled by a factor $K$:*

$$\underline{\Lambda}(n) = \{\underline{x} : B(\underline{x}) \cap \Lambda(n) \neq \varnothing\}.$$

*If there exists $\delta > 0$ such that*

$$\forall \underline{x} \in \mathbb{Z}^2 \qquad P[X(\underline{x}) = 0 | X(\underline{z}), |\underline{x} - \underline{z}|_\infty \geq 3] \leq \delta,$$

*then for any $n, K, \varepsilon$ satisfying $n \geq 6K, \delta < \varepsilon$, we have*

$$P\left[\frac{1}{|\underline{\Lambda}(n)|} \sum_{\underline{x} \in \underline{\Lambda}(n)} 1_{X(\underline{x})=0} \geq \varepsilon\right] \leq 9 \exp\left(-\Lambda^*(\varepsilon, \delta)\left\lfloor \frac{n}{6K}\right\rfloor^2\right),$$

*where*

$$\Lambda^*(\varepsilon, \delta) = \varepsilon \log \frac{\varepsilon}{\delta} + (1 - \varepsilon) \log \frac{1 - \varepsilon}{1 - \delta}$$

*is the Cramér transform of a Bernoulli variable with parameter $\delta$.*

We finish with Hoeffding's inequality that will be useful.

Lemma 7 (Theorem 1 of [23]). *If $(X_i)_{1 \leq i \leq n}$ are independent random variables with values in $[-1, 1]$ and with mean $m$, then*

$$\forall t \in ]0, 1 - m[ \qquad P\left[\sum_{i=1}^n (X_i - m) \geq nt\right] \leq \exp(-nt^2).$$

4.2. *Decoupling and preliminary block estimates near criticality.* This section is devoted to preliminary FK-percolation estimates near criticality. These results are the subject of an independent paper [11] and are included here for self-consistency.

4.2.1. *Decoupling.* In order to decouple distant events near criticality, we rely on an adaptation of weak mixing results in a situation where $p \to p_c$. This adaptation is contained in [11] from which we derive the following decoupling lemma.

LEMMA 8. *Let $p \neq p_c$, $a > 5$, $\Lambda$ a box and $\Phi \in \mathcal{FK}(p, \Lambda)$. There exist two positive constants $\lambda$ and $c = c(a)$ such that for every two sets $\Gamma, \Delta \subset \Lambda$ satisfying*

$$d(\Gamma, \Delta) > c \left( |p - p_c|^{-a} \vee \frac{\log |\Gamma|}{|p - p_c|} \vee \frac{\log |\Delta|}{|p - p_c|} \right)$$

*and for every two events $A \in \mathcal{F}_\Gamma$ and $B \in \mathcal{F}_\Delta$, we have*

$$|\Phi[A \cap B] - \Phi[A]\Phi[B]| \leq \exp[-\lambda(p - p_c)d(\Gamma, \Delta)]\Phi[A]\Phi[B].$$

PROOF. Let $a > 5$. First one needs to have a control of the exponential decay of connectivities in finite boxes when $p \uparrow p_c$ and the size $n$ of the box goes to infinity faster than $(p_c - p)^{-a}$. This is the subject of Proposition 9 of [11]. Once this has been established, an analogue of the weak mixing result contained in Theorem 3.1 of [1] can be obtained. From there, an adaptation of the arguments in Lemma 3.2 of [3] establishes the result. □

4.2.2. *Block estimates near criticality.* This subsection contains the preliminary block estimates established in [11] which are needed to implement a proper coarse graining. In the following we will use the boxes

$$\widetilde{\Lambda}(n) = \{x \in \mathbb{Z}^2 : d(x, \Lambda(n)) \leq n/10\}.$$

We take these bigger boxes in order to give estimates on events that occur in $\Lambda(n)$ uniformly over the boundary conditions on $\widetilde{\Lambda}(n)$. In fact, an adaptation of [2] to our regime would spare us this precaution.

We will say that a FK-cluster $C$ of a box $\Lambda$ is a *crossing* cluster or that $C$ *crosses* the box $\Lambda$ if $C$ connects all the sides of $\Lambda$. Note that in dimension two if there exists a crossing cluster in a box then it is necessarily unique. We will give estimates on the following block events:

$$U(\Lambda) = \{\exists \text{ an open crossing cluster } C^* \text{ in } \Lambda\}.$$

For $M > 0$, we define

$$R(\Lambda, M) = U(\Lambda) \cap \{\text{every open path } \gamma \subset \Lambda \text{ with } \operatorname{diam}(\gamma) \geq M \text{ is in } C^*\}$$

$$\cap \{C^* \text{ crosses every sub-box of } \Lambda \text{ with diameter } \geq M\}.$$

For $\delta > 0$, we define

$$V(\Lambda, \delta) = U(\Lambda) \cap \{|C^*| \geq (1-\delta)\theta|\Lambda|\},$$

$$F(\Lambda, \delta) = \left\{ \begin{array}{l} \exists \text{ an open circuit } \gamma \text{ enclosing a volume } \geq (1-\delta)|\Lambda| \\ \text{and such that } \sup_{x \in \gamma} d(x, \partial\Lambda) \leq \delta|\partial\Lambda| \end{array} \right\},$$

$$W(\Lambda, \delta) = \{|\{x \in \Lambda : x \leftrightarrow \partial\Lambda\}| \leq (1+\delta)\theta|\Lambda|\},$$

$$T(\Lambda, \delta) = \left\{ \left| \sum_{x \in \Lambda : x \leftrightarrow \partial\Lambda} \sigma(x) \right| \leq \delta\theta|\Lambda| \right\}.$$

Notice that the last event involves the FK-Ising coupling. Let us begin with the first two events.

Lemma 9.    *Let $a > 5$. There exist two positive constants $\lambda, c = c(a)$ such that if $p > p_c$ and $n > c(p - p_c)^{-a}$ then*

$$\forall \Phi \in \mathcal{FK}(\widetilde{\Lambda}(n), p) \qquad \log \Phi[U(\Lambda(n))^c] \leq -\lambda(p - p_c)n.$$

*Moreover, if $M$ is such that*

$$(6) \qquad\qquad \frac{\log n}{\kappa(p - p_c)} < M \leq n$$

*with $\kappa > 0$ small enough, then*

$$\forall \Phi \in \mathcal{FK}(\widetilde{\Lambda}(n), p) \qquad \log \Phi[R(\Lambda(n), M)^c] \leq -\lambda(p - p_c)M.$$

Now we turn to the following estimation of the crossing cluster's size:

Lemma 10.    *Let $p > p_c$ and $\delta > 0$. Let $a > 5$ and $\alpha \in ]0, (1 + \frac{1}{8a})^{-1}[$. There exists a positive constant $c = c(a, \alpha)$ such that, if $n \uparrow \infty$ and $p \downarrow p_c$ in such a way that $n^\alpha(p - p_c)^a > c$, then*

$$\sup_{\Phi \in \mathcal{FK}(\widetilde{\Lambda}(n), p)} \Phi[V(\Lambda(n), \delta)^c]$$

$$\leq \exp(-\lambda\delta(p - p_c)n^\alpha) + \exp\left(-\frac{\delta^2\theta^2(p)}{4}n^{2-2\alpha}\right),$$

*where $\lambda$ is a positive constant. In particular,*

$$\lim_{n, p} \inf_{\Phi \in \mathcal{FK}(\widetilde{\Lambda}(n), p)} \Phi[V(\Lambda(n), \delta)] = 1.$$

Next, we consider the deviations from above for the size of the crossing cluster.

LEMMA 11.   *Let $p > p_c$ and $\delta > 0$. If $n > 8 m_{\sup}(\delta, p)/\delta$ where*

$$m_{\sup}(\delta, p) = \inf\left\{ m \geq 1 : \forall n \geq m \ \frac{1}{|\Lambda(n)|} \Phi^{w,p}_{\Lambda(n)}[M_{\Lambda(n)}] \leq (1 + \delta/2)\theta \right\},$$

*then*

$$(7) \qquad \log \Phi^{w,p}_{\Lambda(n)}[W(\Lambda(n), \delta)^c] \leq -\left(\frac{\delta\theta n}{4 m_{\sup}(\delta, p)}\right)^2.$$

*In particular, for every $a > 5/4$, there exists a positive constant $c = c(a, \delta)$ such that whenever $n \uparrow \infty$ and $p \downarrow p_c$ in such a way that $n > c(p - p_c)^{-a}$ then*

$$\limsup_{(n,p)} \frac{1}{(p - p_c)^{2a + 1/4} n^2} \log \Phi^{w,p}_{\Lambda(n)}[W(\Lambda(n), \delta)^c] < 0.$$

LEMMA 12.   *Let $a > 5$ and $\delta > 0$. There exist two positive constants $\lambda, c = c(a, \delta)$ such that for all $p > p_c$ and $n > 4$ such that $n > c(p - p_c)^{-a}$, we have*

$$\log \Phi[F(\Lambda(n), \delta)^c] \leq -\lambda\delta n(p - p_c)$$

*uniformly in $\Phi \in \mathcal{FK}(\widetilde{\Lambda}(n), p)$.*

PROOF.   It suffices to note that

$$\Phi[F(\Lambda(n), \delta)^c] \leq \Phi[\partial\Lambda(n(1 - 2\delta)) \leftrightarrow \partial\Lambda(n) \text{ by an open dual path}]. \quad \square$$

Finally, we consider the event $T$ where the edge and the spin configuration of the FK-Ising coupling $\mathbb{P}^+_n$, defined in Section 3.2, is involved.

LEMMA 13.   *Let $\delta > 0$ and $a > 5$. If $p \downarrow p_c$ and $n \uparrow \infty$ in such a way that $n > (p - p_c)^{-a}$, then*

$$\lim_{n,p} \mathbb{P}^+_n[T(\Lambda(n), \delta)] = 1.$$

PROOF.   Let $\mathcal{C}$ be the collection of the open clusters which do not touch the boundary $\partial\Lambda(n)$. Let $p, n, M = (\log n)/\kappa(p - p_c)$ where $\kappa$ is as in Lemma 9. Let $\omega \in R(\Lambda(n), M)$. Using Chebyshev's inequality,

$$\mathbb{P}^+_n[T(\Lambda(n), \delta)^c | \omega] \leq \mathbb{P}^+_n\left[ \frac{1}{|\Lambda(n)|} \left| \sum_{C \in \mathcal{C}} \sigma(C)|C| \right| \geq \delta\theta \Big| \omega \right]$$

$$\leq \frac{1}{\delta^2\theta^2 |\Lambda(n)|^2} \sum_{C \in \mathcal{C}} |C|^2 \leq \frac{M^4}{\delta^2\theta^2 n^2}.$$

Imposing $\log n/((p - p_c)\sqrt{\theta n}) \to 0$, using Lemma 9 and the previous inequalities, we get the desired result.   $\square$

**5. The upper bound.**     Let $x$ be a point of $\mathbb{R}^2$. The closed ball of center $x$ and Euclidean radius $r > 0$ is denoted by $B(x, r)$. For $w$ in the unit sphere $S^1$, we define the half balls

$$B_-(x, r, w) = B(x, r) \cap \{y \in \mathbb{R}^2 : (y - x) \cdot w \leq 0\},$$

$$B_+(x, r, w) = B(x, r) \cap \{y \in \mathbb{R}^2 : (y - x) \cdot w \geq 0\}.$$

To prove the local upper bound we need to estimate an FK-percolation event which occurs when the locally averaged magnetization exhibits a jump. We will do this by showing that this event implies the existence of an interface. The relevant event is that there exists a collection $\mathcal{G}$ of open clusters in $B(nx, nr)$ such that

$$\sum_{C \in \mathcal{G}} |C \cap B_-(nx, nr, w)| \geq (1 - \theta\delta)\mathcal{L}^2(B_-(nx, nr, w)),$$

$$\sum_{C \in \mathcal{G}} |C \cap B_+(nx, nr, w)| \leq \delta\theta\mathcal{L}^2(B_+(nx, nr, w)).$$

We will denote this event by $\text{Sep}(n, x, r, w, \delta)$. Next, we state the so-called interface lemma whose proof is given after some preliminary work.

LEMMA 14.     *Let $x \in Q$ and $0 < r \leq 1$ such that $B(x, r) \subset Q$. Let $\delta > 0$ and $w \in S^1$. If $p \downarrow p_c$ and $n \uparrow \infty$ in such a way that $n(p - p_c)^{20} \to \infty$ and such that $\log(n)/\log(1/(p - p_c))$ stays bounded, then*

$$\limsup_{(n, p)} \frac{1}{(p - p_c)n} \log \Phi_{\Lambda(n)}^{w, p}[\text{Sep}(n, x, r, w, \delta)]$$

$$\leq -2r\tau_c(1 - c'\delta^{1/2}),$$

*where $c'$ is a positive constant.*

In [9, 10, 12, 13] a cutting procedure has been used to create an interface from the event Sep without altering the probability too much. In our context such an approach does not work. This stems from the fact that the monotone perturbation lemma (Lemma 6.3 of [12]) is not appropriate when $p \downarrow p_c$. We thus have to proceed differently. We start by showing that the event Sep is well approximated by a similar event involving filled clusters instead of clusters with a lot of small holes; then instead of cutting some edges in order to create the interface, we will detect a piecewise interface. Let us fix a small positive $\eta$ that will be determined later and $\rho$ such that

$$0 < \eta < \rho < r, \qquad 0 < 2\eta < \sqrt{r^2 - \rho^2},$$

and we restrict our attention to the rectangle,

$$R = \{y \in B(nx, nr) : -\eta n \leq (y - nx) \cdot w \leq \eta n, -\rho n \leq (y - nx) \cdot w_\perp \leq \rho n\},$$

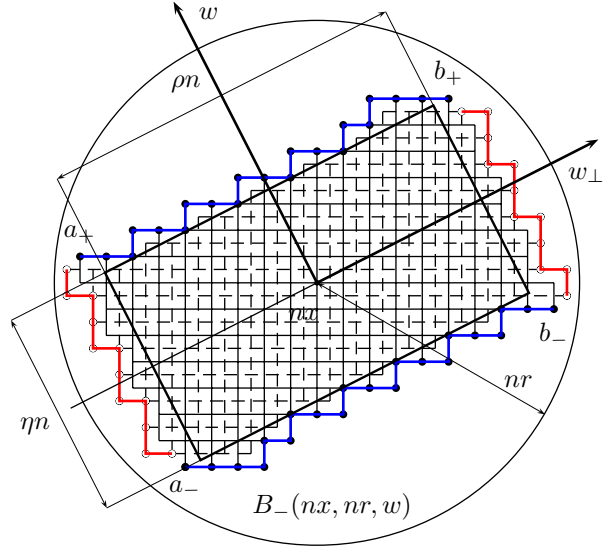

FIG. 1. *The interface.*

where $w_\perp$ is the vector perpendicular to $w$ such that $(w_\perp, w)$ is a direct basis.

We define the *right* to be the direction at which $w_\perp$ points and the *top* the direction at which $w$ points and accordingly we define the *left* and the *bottom*. We consider the graph $D \subset \Lambda(n)$ whose edges are the boundaries of the squares centered in $\mathbb{Z}^2 + (1/2, 1/2)$ that intersect $R$. In this way, the set $\partial D$ is a simple closed circuit. We denote by $a_+$ the upper left site of the square that contains the upper left point of $\partial D$. Going clockwise we define successively and in a similar way $b_+$, $b_-$ and $a_-$, see the Figure 1. We define the top boundary $\partial^+ D$ as the path of $\partial D$ that joins clockwise $a_+$ to $b_+$. Similarly, the bottom boundary $\partial^- D$ is the path of $\partial D$ that joins counterclockwise $a_-$ to $b_-$. We define also $D_+ = B_+(nx, nr, w) \cap D$ and $D_- = B_-(nx, nr, w) \cap D$. Since the interface is an open *dual* path that goes from the left to the right, we need also to consider the dual graph $\widehat{D}$ of $D$ which is depicted in the previous figure. This permits us to define the right boundary $\partial^R \widehat{D}$ as the piece of the boundary $\partial \widehat{D}$ that joins clockwise the center of the square containing $b_+$ to the center of the square containing $b_-$. Similarly we define the left boundary $\partial^L \widehat{D}$. The definitions of $\partial^+ D, \partial^- D, \partial^L \widehat{D}$ and $\partial^R \widehat{D}$ guarantee that, if a configuration $\omega \in \Omega_D$ does not contain any open cluster that connects $\partial^+ D$ to $\partial^- D$ then, in the dual configuration $\widehat{\omega} \in \Omega_{\widehat{D}}$, there exists an open cluster that connects $\partial^L \widehat{D}$ to $\partial^R \widehat{D}$.

In order to prove the upper bound, it is sufficient to consider the consequence of $\mathrm{Sep}(\delta)$ on the configuration restricted to $D$ which is a convenient set for duality arguments. To depict our restriction to $D$, we denote by $\mathcal{C}$ the clusters in $D$ that connect $\partial^+ D$ to $\partial^- D$ and suppose that there exists a collection $\mathcal{G}$ of open clusters in $B(nx, nr, w)$ that realizes $\mathrm{Sep}(\delta)$. In this situation, we can make the decomposition $\mathcal{C} = \mathcal{C}_- \cup \mathcal{C}_+$ where $\mathcal{C}_-$ is the collection of the open clusters in $\mathcal{C}$ that are contained in a cluster of the collection $\mathcal{G}$, and $\mathcal{C}_+$ is the collection of the open clusters in $\mathcal{C}$ that are not contained in any cluster of the collection $\mathcal{G}$. By the definition of the event Sep, the cardinality of the intersection of $D_+$ with the clusters in $\mathcal{C}_-$ is less than $\theta \delta \mathcal{L}^2(B_+(nx, nr, w))$, and the cardinality of the intersection of $D_-$ with the clusters in $\mathcal{C}_+$ is also less than $\theta \delta \mathcal{L}^2(B_-(nx, nr, w))$. Thus the event Sep implies the following event involving only the clusters in $D$:

$$\mathrm{Sep}_D(\delta) = \Bigg\{ \text{there exists a decomposition } \mathcal{C} = \mathcal{C}_- \cup \mathcal{C}_+ \text{ such that}$$

$$(8) \qquad \sum_{C \in \mathcal{C}_-} |C \cap D_+(nx, nr, w)| \leq \pi \delta \theta (nr)^2,$$

$$\sum_{C \in \mathcal{C}_+} |C \cap D_-(nx, nr, w)| \leq \pi \delta \theta (nr)^2 \Bigg\}.$$

5.1. *Elimination of the small holes.* In our situation it is necessary to fill the small holes of the clusters that enter in the definition of Sep in order to give an adequate upper bound. Namely, we will replace the event Sep with an event Sep′ that uses only filled clusters and we will show that the probability of Sep is well approximated by the probability of the auxiliary event Sep′. In order to construct our filling procedure we need some definitions. Let $C \in \mathcal{C}$. We introduce the notion of *holes* of $C$. For this we consider the dual $\mathbb{E}(\widehat{D}) \setminus \mathbb{E}(C)$ of the complement of $C$. Each maximally connected set $\widehat{F}$ of $\mathbb{E}(\widehat{D}) \setminus \mathbb{E}(C)$ that is isolated from the other clusters of $\mathcal{C}$ by $C$ will be designated as a *hole* of $C$. For each hole we also define the following notion of boundary:

$$\Delta \widehat{F} = \{ \widehat{e} \in \widehat{F} : d(\widehat{e}, \mathbb{E}(C)) = \tfrac{1}{2} \}.$$

Note that by the definition of a hole, the edges of $\Delta \widehat{F}$ are all open dual edges.

Next, we fix $M < n$ and say that a hole $\widehat{F}$ is small (respectively big) if $\mathrm{diam}(\widehat{F}) < M$ [respectively, if $\mathrm{diam}(\widehat{F}) \geq M$]. For each $C \in \mathcal{C}$ we define its filling $\mathrm{fill}\, C$ as

$$\mathrm{fill}\, C = C \cup \bigcup_{\widehat{F}} F,$$

where the union runs over all the small holes of $C$ and where $F$ is the set of edges in $D$ whose dual is $\widehat{F}$. Note that if $C_1 \neq C_2$ then $\mathrm{fill}\, C_1 \neq \mathrm{fill}\, C_2$.

For $* = -, +$, let

$$\mathcal{C}^{\text{fill}}_* = \{\text{fill } C : C \in \mathcal{C}_*\}.$$

We define then a modified $\text{Sep}_D$ event that involves only the following filled clusters:

$$\text{Sep}'_D(n, x, r, w, \delta) = \left\{ \text{there exists a decomposition } \mathcal{C} = \mathcal{C}_- \cup \mathcal{C}_+ \text{ such that} \right.$$

$$\left. \sum_{S \in \mathcal{C}^{\text{filled}}_-} |S \cap D_+| \leq \delta \pi (nr)^2, \sum_{S \in \mathcal{C}^{\text{filled}}_+} |S \cap D_-| \leq \delta \pi (nr)^2 \right\}.$$

Note that the event $\text{Sep}'_D$ involves only the filled clusters of $D$, and even if Sep has been defined originally in $B(nx, nr)$, we will only use its consequence (8). Now we show that the event $\text{Sep}_D$ is well approximated by the event $\text{Sep}'_D$:

LEMMA 15. *Let $\delta > 0$, $a > 5$ and $\alpha \in ]0, (1 + \frac{1}{8a})^{-1}[$. There exists a positive constant $c = c(a, \alpha, \delta) > 0$ such that if*

$$\eta r n > M > c(p - p_c)^{-a/\alpha},$$

*then uniformly in $n, M, x, r, w$ we have that*

$$\log \Phi^{w,p}_{\Lambda(n)}[\text{Sep}_D(n, x, r, w, \delta) \setminus \text{Sep}'_D(n, x, r, w, 4\delta)]$$

$$\leq -\lambda \frac{n^2}{M^2} \delta (p - p_c)^{1/8} \log\left(\frac{1}{p - p_c}\right),$$

*where $\lambda$ is a positive constant.*

PROOF. First we renormalize $D$ into $\underline{D}$ by partitioning it with blocks $B(\underline{x})$ of size $M/2$. We say that a block $B(\underline{x})$ is good if and only if the event

$$V(B(\underline{x}), \delta) \cap W(B(\underline{x}), \delta) \cap F(B(\underline{x}), \delta) \cap R(B'(\underline{x}), M/4)$$

occurs. Recalling the definitions of $V, W, F, R$, the above event is

$$\left\{ \begin{array}{l} \exists \text{ crossing cluster } C^* \text{ in } B(\underline{x}) \\ \text{and } |C^*| \geq (1 - \delta)\theta |B(\underline{x})| \end{array} \right\}$$

$$\cap \{|\{x \in B(\underline{x}) : x \leftrightarrow \partial B(\underline{x})\}| \leq (1 + \delta)\theta |B(\underline{x})|\}$$

$$\cap \left\{ \begin{array}{l} \exists \text{ an open circuit } \gamma \text{ in } B(\underline{x}) \text{ enclosing a volume } \geq (1 - \delta)|B(\underline{x})| \\ \text{and such that } \sup_{x \in \gamma} d(x, \partial B(\underline{x})) \leq \delta |\partial B(\underline{x})| \end{array} \right\}$$

$$\cap \{\exists \text{ crossing cluster } \widetilde{C}^* \text{ in } B'(\underline{x})\}$$

$$\cap \left\{ \begin{array}{l} \text{Every open path } \gamma' \subset B'(\underline{x}) \text{ with} \\ \text{diam } \gamma' \geq M/4 \text{ is included in } \widetilde{C}^* \end{array} \right\}$$

$$\cap \{\widetilde{C}^* \text{ crosses every sub-box of } B'(\underline{x}) \text{ of diam} \geq M/4\}.$$

We define the block process $(X(\underline{x}), \underline{x} \in \underline{D})$ by $X(\underline{x}) = 1_{B(\underline{x}) \text{ is good}}$ for $\underline{x} \in \underline{D}$. We have

$$\sum_{C \in \mathcal{C}_-} |C \cap D_+| \leq \sum_{\underline{x} \in \underline{D}_+} \sum_{C \in \mathcal{C}_-} |C \cap B(\underline{x})|$$

$$= \sum_{\substack{\underline{x} \in \underline{D}_+ \\ B(\underline{x}) \text{ is good}}} \sum_{C \in \mathcal{C}_-} |C \cap B(\underline{x})| + \sum_{\substack{\underline{x} \in \underline{D}_+ \\ B(\underline{x}) \text{ is bad}}} \sum_{C \in \mathcal{C}_-} |C \cap B(\underline{x})|.$$

When $B(\underline{x})$ is good and $\mathrm{fill}(C) \cap B(\underline{x}) \neq \varnothing$ then $C \cap B(\underline{x}) \neq \varnothing$. We also have that

$$(9) \qquad |C \cap B(\underline{x})| \geq |C^*| \geq (1 - \delta)\theta|B(\underline{x})|.$$

The first inequality holds because when $B(\underline{x})$ is good, then every cluster $C \in \mathcal{C}_-$ that intersects $B(\underline{x})$ will contain a path of diameter at least $M/2$ in the box $B'(\underline{x})$; this path is included in $\widetilde{C}^*$, and thus the cluster $C$ contains $\widetilde{C}^*$ which also contains $C^*$. Next

$$(10) \qquad |C \cap B(\underline{x})| \leq |\{x \in B(\underline{x}) : x \leftrightarrow \partial B(\underline{x})\}| \leq (1 + \delta)\theta|B(\underline{x})|.$$

This inequality is true because for reasons of diameter no cluster $C \in \mathcal{C}_-$ fits into a box $B(\underline{x})$. Thus all the connected components of $C \cap B(\underline{x})$ have to be connected to $\partial B(\underline{x})$. Next, as before, when $B(\underline{x})$ is good then $\widetilde{C}^* \cap B(\underline{x}) \subset C \cap B(\underline{x})$. And since $\mathrm{diam}\,\gamma > M/4$, we have that $\gamma \subset \widetilde{C}^* \cap B(\underline{x}) \subset C \cap B(\underline{x})$. But the diameter of $\gamma$ is less than $M$, and thus the interior of $\gamma$ is included in $\mathrm{fill}\,C$; hence

$$(11) \qquad (1 - \delta)\theta|B(\underline{x})| \leq \theta|\mathrm{fill}(C) \cap B(\underline{x})| \leq \theta|B(\underline{x})|.$$

By (9), (10), (11) we get

$$-2\delta\theta|B(\underline{x})| \leq |\mathrm{fill}(C) \cap B(\underline{x})|\theta - |C \cap B(\underline{x})| \leq 2\delta\theta|B(\underline{x})|.$$

Since two different clusters of $\mathcal{C}$ cannot intersect the same good block, we obtain

$$\sum_{C \in \mathcal{C}_-} ||C \cap D_+| - \theta|\mathrm{fill}(C) \cap D_+|| \leq \sum_{\underline{x} \text{ good}} 2\delta\theta|B(\underline{x})| + \sum_{\underline{x} \text{ bad}} |B(\underline{x})|$$

$$\leq 2\delta\theta\pi(nr)^2 + \frac{|\{\underline{x} \in \underline{D} : X(\underline{x}) = 0\}|}{|\underline{D}|} 2\pi(nr)^2.$$

Doing the same reasoning for $D_-$ with $\mathcal{C}_+$, we get

$$\sum_{C \in \mathcal{C}_+} ||C \cap D_-| - \theta|\mathrm{fill}(C) \cap D_-|| \leq 2\delta\theta\pi(nr)^2 + \frac{|\{\underline{x} \in \underline{D} : X(\underline{x}) = 0\}|}{|\underline{D}|} 2\pi(nr)^2.$$

From this, we conclude that

$$\mathrm{Sep}_D(\delta) \cap \left\{ \frac{|\{\underline{x} \in \underline{D} : X(\underline{x}) = 0\}|}{|\underline{D}|} \leq \frac{\delta\theta}{2} \right\}$$

is included in $\mathrm{Sep}'_D(4\delta)$. Thus

$$\Phi^{w,p}_{\Lambda(n)}[\mathrm{Sep}_D(n,x,r,w,\delta) \setminus \mathrm{Sep}'_D(n,x,r,w,4\delta)]$$

$$\leq \Phi^{w,p}_{\Lambda(n)}\left[ \frac{|\{\underline{x} \in \underline{D} : X(\underline{x}) = 0\}|}{|\underline{D}|} \geq \frac{\delta\theta}{2} \right].$$

Finally, we show that it is possible to tune our regime so that with probability very close to one, the fraction of bad boxes in $D$ remains negligible. Fix $a > 5$ and $\alpha \in ]0, (1 + \frac{a}{8})^{-1}[$. There exists a positive $c = c(\alpha, a, \delta)$ such that by Lemmas 9, 10, 11, 12,

$$\sup_{\Phi \in \mathcal{FK}(p,D)} \Phi[X(\underline{x}) = 0 | X(\underline{y}), |\underline{x} - \underline{y}|_\infty \geq 3] = \rho(M, p) \downarrow 0,$$

when $M \uparrow \infty$ and $p \downarrow p_c$ in such a way that $M^\alpha > c(p - p_c)^{-a}$. Thus, by Lemma 6 we get that

$$\log \Phi^{w,p}_{\Lambda(n)}[\mathrm{Sep}_D(n,x,r,w,\delta) \setminus \mathrm{Sep}'_D(n,x,r,w,4\delta)] \leq -\delta\theta \log \frac{\delta\theta}{\rho(M,p)} \left\lfloor \frac{nr}{3M} \right\rfloor^2.$$

By using Onsager's formula we get that $\theta \sim (p - p_c)^{1/8}$ when $p \downarrow p_c$. The conclusion follows from the speed of the convergence $\rho(M,p) \downarrow 0$ provided by Lemmas 9, 10, 11, 12.   $\square$

5.2. *The piecewise interface.* In this section we will detect a piecewise interface from the occurrence of the event $\mathrm{Sep}'_D$. We suppose that the event $\mathrm{Sep}'_D$ occurs, and let $\mathcal{C}_-, \mathcal{C}_+$ be a decomposition of $\mathcal{C}$ realizing it. We define

$$\mathcal{C}^{\mathrm{filled}} = \mathcal{C}^{\mathrm{filled}}_- \cup \mathcal{C}^{\mathrm{filled}}_+,$$

where $\mathcal{C}^{\mathrm{filled}}_*$ has been defined in Section 5.1. Note that there is a natural order in $\mathcal{C}^{\mathrm{filled}}$ and thus it is possible to enumerate the elements of $\mathcal{C}^{\mathrm{filled}}$ from the left to the right from $S_1$ to $S_{|\mathcal{C}^{\mathrm{filled}}|}$. Next, for each $h \in \mathbb{R}$, we define the line

$$\pi(h) = \{y \in \mathbb{R}^2 : (y - nx) \cdot w = h\}.$$

Let

$$E_+ = \bigcup_{S \in \mathcal{C}^{\mathrm{filled}}_-} \mathbb{E}(S \cap D_+)$$

be the set of the edges in $D_+$ that belong to a filled cluster of $\mathcal{C}_-^{\text{filled}}$. We define similarly $E_-$ as the set of the edges in $D_-$ that belong to a filled cluster of $\mathcal{C}_+^{\text{filled}}$. Then

$$\int_{\eta n/3}^{2\eta n/3} |\{e \in E_+ : e \cap \pi(h) \neq \varnothing\}| \, dh \leq |E_+|,$$

where $\eta$ has been defined in the paragraph after Lemma 14. Since $\mathcal{C}_-$, $\mathcal{C}_+$ realize the event $\text{Sep}'_D(n, x, r, w, \delta)$, then we have that $|E_+| \leq \delta\pi(nr)^2$ and by the previous inequality, there exists $h \in [\eta n/3, 2\eta n/3]$ such that

$$(12) \qquad |\{e \in E_+ : e \cap \pi(h) \neq \varnothing\}| \leq \frac{3\delta}{\eta} n\pi r^2.$$

Let $h^*$ be the infimum in $[\eta n/3, 2\eta n/3]$ of the real numbers $h$ satisfying this inequality. If we increase the value of $h$ by a small $\varepsilon > 0$, then the inequality (12) still holds and $\pi(h^* + \varepsilon) \cap \mathbb{Z}^2 \cap D = \varnothing$. We choose one such $h^* + \varepsilon$ and we call it $h_+$. Moreover any edge of $E_+$ which intersects $\pi(h_+)$ has an endpoint in each of the two half spaces delimited by $\pi(h_+)$. In a symmetric way we get from $E_-$ a value $h_-$ in $[-2\eta n/3, -\eta n/3]$. The edges in $\{e \in E_- : e \cap \pi(h_-) \neq \varnothing\} \cup \{e \in E_+ : e \cap \pi(h_+) \neq \varnothing\}$ will be designated as *bad edges*. We end with a horizontal segment $\pi(h_+)$ in $D_+$ that crosses at most $3\pi\delta n r^2/\eta$ edges belonging to a cluster $C \in \mathcal{C}_-$ and a horizontal segment $\pi(h_-)$ in $D_-$ that crosses at most $3\pi\delta n r^2/\eta$ edges belonging to a cluster $C \in \mathcal{C}_+$. Note that if these $6\pi\delta n r^2/\eta$ bad edges were closed then by duality, there would exist an open dual path connecting $\partial^L \widehat{D}$ to $\partial^R \widehat{D}$. For $* = -, +$, we introduce the following sets of edges:

$$\Pi_* = \{e \in \mathbb{E}(D) : e \cap \pi(h_*) \neq \varnothing\}, \qquad \widehat{\Pi}_* = \{\widehat{e} \in \mathbb{E}(\widehat{D}) : e \in \Pi_*\}.$$

The set $\Pi_*$ is the set of all the edges that intersect $\pi(h_*)$ and $\widehat{\Pi}_*$ is its dual set. Note that $\widehat{\Pi}_*$ is always a simple dual path connecting $\partial^L \widehat{D}$ to $\partial^R \widehat{D}$.

In order to capture the relevant dual connections, we introduce for each dual path $\gamma \subset \widehat{D}$, its $w$-diameter:

$$\text{diam}_w(\gamma) = \max_{x, y \in \gamma} (y - x) \cdot w_\perp.$$

Lemma 16. *If the event* $\text{Sep}'_D(n, x, r, w, \delta)$ *occurs, then there exists a family of open dual paths* $(\widehat{\xi}_1, \widehat{\xi}_2, \ldots, \widehat{\xi}_K)$ *such that*

$$\text{diam}_w(\widehat{\xi}_1) + \cdots + \text{diam}_w(\widehat{\xi}_K) \geq 2n\rho - \frac{6\delta}{\eta} n\pi r^2,$$

*where $\rho$ has been defined just after Lemma 14. Moreover the number $K - 1$ is bounded above by the number of open dual clusters of diameter $\geq M$ that intersect $\pi(h_+) \cup \pi(h_-)$.*

PROOF.   First choose $\omega \in \mathrm{Sep}'(n, x, r, w, \delta)$. If there is no top–bottom crossing cluster in $D$ then by duality, there exists an open dual path connecting $\partial^L \widehat{D}$ to $\partial^R \widehat{D}$, and we are done. If there are crossing clusters, then we construct an algorithm that detects in every dual configuration $\widehat{\omega}$ of $\omega \in \mathrm{Sep}'(n, x, r, w, \delta)$, a way to move from $\partial^L \widehat{D}$ to $\partial^R \widehat{D}$ by using either open dual paths or paths of bad edges, that is, edges of $(\widehat{\Pi}_+ \cap \widehat{E}_+) \cup (\widehat{\Pi}_- \cap \widehat{E}_-)$. Note that the paths of bad edges are not necessarily open dual paths and we will denote them by *tunnels*. Using tunnels will mean following a path of bad edges along $\widehat{\Pi}_+ \cup \widehat{\Pi}_-$ from the left to the right until we reach an edge that is not bad. Note already that the total length of the tunnels is bounded by

$$|(\widehat{\Pi}_+ \cap \widehat{E}_+) \cup (\widehat{\Pi}_- \cap \widehat{E}_-)| \le 6\delta\pi r^2 n/\eta.$$

Let us first sketch the idea behind the algorithm: we want to move from the left to the right along open dual paths. The only obstacles preventing us from doing that is the existence of top–bottom crossing clusters. To overcome the problem, when we meet such a filled crossing cluster $S$, we check if $S$ is in $\mathcal{C}_-^{\mathrm{fill}}$ or in $\mathcal{C}_+^{\mathrm{fill}}$. Accordingly, we move to $\widehat{\Pi}_+$ or $\widehat{\Pi}_-$ and traverse the obstacle using a tunnel. After such a tunnel, we meet holes $\widehat{F}$ of $S$ that we traverse using open dual paths included in $\Delta\widehat{F}$. We continue like this until we reach the right-hand side of $S$. After this, we find an open dual path that reaches the next top–bottom crossing cluster and so on. At the end, the total number of closed dual edges that we have used is negligible. But this is not enough as the number of segments of open dual paths may be very large and this may prevent us from decoupling properly the probability of these segments. It is at this point that the filling of the small holes is important. Indeed, with our filling, we are guaranteed that at each time we produce a new open dual path, we will meet a large open dual cluster that intersects $\widehat{\Pi}_+ \cap \widehat{\Pi}_-$. The number of such clusters can be controlled in order to decouple the relevant dual connections.

Next, we give the precise description of our algorithm:

*Initialization.*   First we check the leftmost edge $e_{1+}$ of $\Pi_+$ and the leftmost edge $e_{1-}$ of $\Pi_-$.

(1) If $e_{1+}$ is in $\mathcal{C}_-^{\mathrm{fill}}$, then we use the tunnel included in $\widehat{\Pi}_+$ that starts in $\widehat{e}_{1+}$ and ends at an edge of $\widehat{\Pi}_+$ that is not bad.

(2) If $e_{1+}$ is in $\mathcal{C}_+^{\mathrm{fill}}$, then two subcases arise according to $e_{1-}$:

    (2a) If $e_{1-} \in S_1$, where $S_1$ is the first, from the left, top–bottom filled crossing cluster. We use the tunnel included in $\widehat{\Pi}_-$ that starts at $\widehat{e}_{1-}$ and ends at an edge which is not bad.

(2b) If $e_{1-} \notin S_1$, then $e_{1-}$ is isolated from $\partial^R \widehat{D}$ by $S_1$. So there exists an open dual path from $\partial^L \widehat{D}$ to a site in $\widehat{\Pi}_-$. Let $\widehat{\xi}_1$ be such a path whose endpoint on $\widehat{\Pi}_-$ is rightmost. By doing so, the right successor edge of $\widehat{\xi}_1$ on $\widehat{\Pi}_-$ must be a bad edge and thus the entrance of a tunnel. We use this tunnel until we reach an edge that is not bad.

*Intermediate steps.* Next, we suppose that we have reached an edge $\widehat{e}_j \in \widehat{\Pi}_+ \cup \widehat{\Pi}_-$ that is not a bad edge and describe how to proceed with the algorithm in order to reach an other edge $\widehat{e}_{j+1} \in \widehat{\Pi}_+ \cup \widehat{\Pi}_-$ that is not bad. Whether $\widehat{e}_j$ is in $\widehat{\Pi}_+$ or $\widehat{\Pi}_-$ is completely symmetric. We thus suppose that $\widehat{e}_j \in \widehat{\Pi}_+$ and the other case can be deduced by symmetry. If $\widehat{e}_j$ is described as above, then it is either included in a hole of a cluster in $\mathcal{C}^{\text{fill}}$ or it is at the right of the rightmost top–bottom crossing path of a cluster in $\mathcal{C}^{\text{fill}}$.

(1) If $\widehat{e}_j$ is in a hole $\widehat{F}$ of a filled cluster $S \in \mathcal{C}^{\text{fill}}$, then we choose a path from $\Delta \widehat{F}$ that takes us to the rightmost edge of $\widehat{\Pi}_+$, and then we are again at the entrance of a tunnel that we cross. We denote by $\widehat{e}_{j+1}$ the successor on $\widehat{\Pi}_+$ of the exit of the just traversed tunnel. By definition $\widehat{e}_{j+1}$ is not a bad edge.

(2) $\widehat{e}_j$ is just at the right of the rightmost top–bottom crossing path of a filled cluster $S \in \mathcal{C}^{\text{fill}}$. Let $S'$ be the next filled crossing cluster of $\mathcal{C}^{\text{fill}}$. If $S' \in \mathcal{C}^{\text{fill}}_+$, then we go along an open dual path that joins $\widehat{e}_j$ to the rightmost dual intersection $\widehat{e}'_{j+1}$ of $\widehat{\Pi}_-$ with the top–bottom crossing open dual path that is just on the left of the leftmost top–bottom crossing path of $S'$. The edge $\widehat{e}'_{j+1}$ is the entrance of a tunnel that we take until we reach an edge in $\widehat{\Pi}_-$ that is not a bad edge. We call this edge $\widehat{e}_{j+1}$.

*The final step.* It is reached when an edge of $\partial^R \widehat{D}$ has been seen. This must happen in a finite number of steps since we explore partially without repetition the edges of $\widehat{\Pi}_- \cup \widehat{\Pi}_+$ from the left to the right. The number of the edges in the tunnels is bounded by $6\delta\pi r^2 n/\eta$; thus the created open dual paths $\widehat{\xi}_1, \ldots, \widehat{\xi}_K$ satisfy

$$\text{diam}_w(\widehat{\xi}_1) + \cdots + \text{diam}_w(\widehat{\xi}_K) \geq 2n\rho - \frac{6\delta}{\eta} n\pi r^2.$$

In addition to that, the just described algorithm has the property that the creation of a new open dual path corresponds to an additional open dual cluster of diameter larger than $M$ that intersects $\widehat{\Pi}_+ \cup \widehat{\Pi}_-$. Thus we can bound $K - 1$ as stated in the lemma. $\square$

5.3. *Separating the pieces of the interface.* In order to get the right probabilistic upper bound from the existence of the piecewise interface, we have to factorize the probability of the dual connections obtained in Lemma 16 without altering our estimates too much. If we were working in independent Bernoulli percolation then we would simply apply the van den Berg–Kesten inequality. Unfortunately this inequality does not hold in dependent FK-percolation models. To decouple our events, we start by constructing a new family of paths from $(\widehat{\xi}_1, \ldots, \widehat{\xi}_K)$. The new paths will be well separated from each other by a distance of at least $0 < \ell < \delta n$. In order to simplify the notation, we consider without loss of generality that our domain $\widehat{D}$ is centered at the origin. For $-\rho n < h < \rho n$ we define the line $v(h)$ parallel to $w$ and at a relative distance $hn\rho$ from the origin,

$$v(h) = \{x \in \mathbb{R}^2 : x \cdot w_\perp = h\}.$$

In addition to that we will need vertical strips that separate the events, so we define for every site $x \in \mathbb{R}^2$ the strip of width $\ell$ on the right of $x$,

$$H_\ell(x) = \{y \in \mathbb{R}^2 : 0 < (y - x) \cdot w_\perp < \ell\}.$$

Now we give the construction of our new well-separated dual open paths. First we start with the value $h = -n$ and we increase $h$ until the first time we find at least one dual open path $\widehat{\gamma}$ that satisfies:

  (i) $\widehat{\gamma}$ is part of one of the paths of the piecewise interface $(\widehat{\xi}_1, \ldots, \widehat{\xi}_K)$;

  (ii) $\widehat{\gamma}$ starts at a site on $v(h)$ and does not intersect the left half plane defined by $v(h)$;

  (iii) $\operatorname{diam}_w(\widehat{\gamma}) \geq \ell$.

Let us call $h_1$ the first value of $h$ where we stopped. Since $\ell < \delta n$, it is clear that $h_1 < n$ as soon as $\operatorname{Sep}'_D$ occurs. Let us pick, among the above mentioned paths, a path $\widehat{\gamma}$ of maximal $w$-diameter. On $\widehat{\gamma}$ we choose two sites $\widehat{x}_1, \widehat{y}_1 \in \widehat{\gamma}$ that satisfy

$$(\widehat{y}_1 - \widehat{x}_1) \cdot w_\perp = \operatorname{diam}_w(\widehat{\gamma})$$

and we define $\widehat{\gamma}_1$ as a dual open path that joins $\widehat{x}_1$ to $\widehat{y}_1$. Right after this path we put the strip $H_\ell(\widehat{y}_1)$.

Now we suppose that $\widehat{\gamma}_1, \ldots, \widehat{\gamma}_j$ and $H_\ell(\widehat{y}_1), \ldots, H_\ell(\widehat{y}_j)$ have been constructed. Then we start with the value $h = h_j + \ell = y_j \cdot w_\perp + \ell$, we increase $h$ until we find a path $\widehat{\gamma}$ that satisfies the above three criteria, (i)–(iii), and we define $\widehat{\gamma}_{j+1}$ in the same way that we defined $\widehat{\gamma}_1$.

We continue this process until we reach the boundary $\partial^R \widehat{D}$.

After this construction, we end with a sequence of strips separating a family of dual open paths $(\widehat{\gamma}_1, \ldots, \widehat{\gamma}_{K'})$ (see Figure 2).

The constructed paths verify the following:

- For every $1 \leq j \leq K'$, we have that $\mathrm{diam}_w(\widehat{\gamma}_j) \geq \ell$.
- The number $K'$ of the new paths is bounded above by the number of paths $K$ in the original piecewise interface. Indeed, two different paths $\widehat{\gamma}_i, \widehat{\gamma}_j$ cannot be part of the same path $\widehat{\xi}$ of the original interface because when defining the paths $\widehat{\gamma}_j$ we always choose one with maximal $w$-diameter.
- The total $w$-diameter of the new family of paths satisfies

$$\tag{13} \sum_{k=1}^{K'} \mathrm{diam}_w(\widehat{\gamma}_k) \geq 2n\rho - \frac{6\delta}{\eta} n\pi r^2 - 2K\ell.$$

Indeed, we lost from the original total $w$-diameter only for two reasons. The first reason is the fact that we have chosen paths of $w$-diameter larger than $\ell$; this gives a maximal loss of $\ell K$. The second reason is the fact that the construction separates the pieces of the interface by strips. These strips are of width $\ell$, and this gives in the worst case another loss of $\ell K$.

For each $j$, we denote by $\widehat{\Lambda}_j$ the region of $\widehat{D}$ between $H_\ell(y_j)$ and $H_\ell(y_{j+1})$, and for each $k > 0$, we define $\Xi(k)$ as the set of families $(\widehat{\Lambda}_1, \ldots, \widehat{\Lambda}_k)$ that partition the set $\widehat{D}$ as above. Also we define $\Upsilon(k)$ as the set of the families $(s_1, \ldots, s_k) \subset \mathbb{R}^k$ such that

$$\tag{14} \forall j \in \{1, \ldots, k\}, s_j \geq \ell, \qquad \sum_{j=1}^{k} s_j \geq 2n\rho - \frac{6\delta}{\eta} n\pi r^2 - 2k\ell.$$

From Lemma 16 and from the last construction, we get the following result.

LEMMA 17. *Suppose that the event* $\mathrm{Sep}'_D(n, x, r, w, \delta)$ *occurs, and let $K$ be the number of open dual clusters of diameter larger than $M$ that cross* $\pi(h_+) \cup \pi(h_-)$.

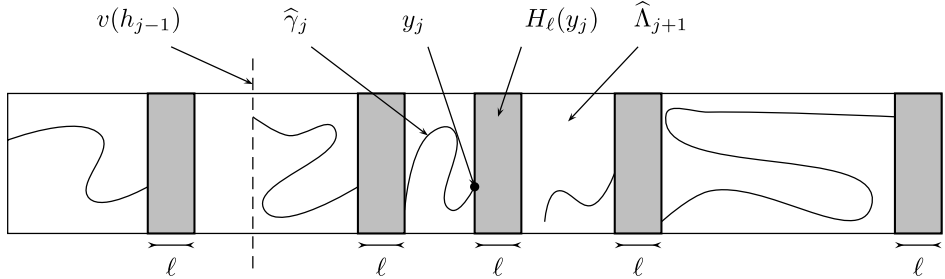

FIG. 2. *The separated pieces of the interface.*

*Then there exist $k \leq K$, $(\widehat{\Lambda}_1, \ldots, \widehat{\Lambda}_k) \in \Xi(k)$ and $(s_1, \ldots, s_k) \in \Upsilon(k)$ such that the event*

$$A(\widehat{\Lambda}_1, \ldots, \widehat{\Lambda}_k; s_1, \ldots, s_k)$$
$$= \bigcap_{j=1}^{k} \{\exists \text{ an open dual path } \widehat{\gamma}_j \subset \widehat{\Lambda}_j \text{ with } \mathrm{diam}_w \, \widehat{\gamma}_j = s_j\}$$

*occurs.*

5.3.1. *Control of the big dual clusters.* Let $h \in (-2r\eta n/3, -r\eta n/3) \cup (r\eta n/3, 2r\eta n/3)$. In what follows, we estimate the number of big open dual clusters that intersect the set $\widehat{\Pi}(h) = \{\widehat{e} \in \mathbb{E}(\widehat{D}) : e \cap \pi(h) \neq \varnothing\}$. Here *big cluster* means a cluster whose diameter exceeds a certain threshold $M > 0$. This estimate is crucial in order to decouple the different pieces of our spatially separated piecewise interface.

LEMMA 18. *Let $p > p_c$ and fix $h \in (-2r\eta n/3, -r\eta n/3) \cup (r\eta n/3, 2r\eta n/3)$. Let $K(h)$ be the number of big open dual clusters of $\widehat{D}$ intersecting $\widehat{\Pi}(h)$. If $\xi > 0$ and $a > 2\xi + 1$ then there exist two positive constants $c = c(a, \xi)$ and $\lambda = \lambda(a, \xi)$ such that*

$$c(p - p_c)^{-a} < M < r\eta n/3$$
$$\Rightarrow \quad \log \Phi[K(h) \geq (p - p_c)^{\xi} n] \leq -\lambda(p - p_c)^{4\xi + 2} nM$$

*uniformly over $\Phi \in \mathcal{FK}(p, D)$.*

PROOF. For a given $h \in (-2r\eta n/3, -r\eta n/3) \cup (r\eta n/3, 2r\eta n/3)$, let

$$\widehat{\Pi}_+(M) = \{y \in \widehat{D} : 0 \leq (y - nx) \cdot w - h \leq M, |(y - nx) \cdot w_\perp| \leq n + M\},$$
$$\widehat{\Pi}_-(M) = \{y \in \widehat{D} : -M \leq (y - nx) \cdot w - h \leq 0, |(y - nx) \cdot w_\perp| \leq n + M\}.$$

Let $\widehat{C}$ be an open dual cluster of diameter $\geq M$ which intersects $\widehat{\Pi}(h)$. Since $|\widehat{C}| \geq M$, either $\widehat{\Pi}_-(M) \cap \widehat{C}$ contains at least $M/2$ sites connected inside $\widehat{\Pi}_-(M)$ to $\partial \widehat{\Pi}_-(M)$ or $\widehat{\Pi}_+(M) \cap \widehat{C}$ contains at least $M/2$ sites connected inside $\widehat{\Pi}_+(M)$ to $\partial \widehat{\Pi}_+(M)$. Thus for $\xi > 0$, we have that

$$\Phi[K \geq (p - p_c)^{\xi} n] \leq \Phi\left[\widehat{M}_{\widehat{\Pi}_-(M)} \geq (p - p_c)^{\xi} \frac{Mn}{4}\right]$$
$$+ \Phi\left[\widehat{M}_{\widehat{\Pi}_+(M)} \geq (p - p_c)^{\xi} \frac{Mn}{4}\right],$$

where $\widehat{M}_{\widehat{\Lambda}} = |\{\widehat{x} \in \widehat{\Lambda} : \widehat{x} \leftrightarrow \partial\widehat{\Lambda}\}|$. We choose $a > 2\xi + 1$ and partition $\widehat{\Pi}_-(M)$ and $\widehat{\Pi}_+(M)$ into blocks $B(\underline{x})$ of size $m = c(p - p_c)^{\xi - a}/4$. We assume that

$$(15) \qquad \min\left(\frac{1}{6}, \frac{\eta r}{3}\right) n > M > 10c(p - p_c)^{-a},$$

where $c$ is a positive constant that will be determined later. Next, we define for $* = +, -,$

$$\widehat{\Pi}'_*(M) = \bigcup_{\underline{x} \in \mathbb{Z}^2 \, : \, B(\underline{x}) \cap \widehat{\Pi}_*(M) \neq \varnothing} B(\underline{x}).$$

Also, the number of partitioning blocks $|\underline{\widehat{\Pi}}'_*(M)|$ satisfies

$$\frac{nM}{m^2} \leq |\underline{\widehat{\Pi}}'_*(M)| \leq \frac{2(M + m)(n + M + m)}{m^2} \leq 7\frac{nM}{m^2}.$$

By subadditivity, one gets

$$\frac{\widehat{M}_{\widehat{\Pi}_*(M)}}{Mn} \leq \frac{7}{|\underline{\widehat{\Pi}}'_*(M)|} \sum_{\underline{x} \in \widehat{\underline{\Pi}}'_*(M)} \frac{\widehat{M}_{B(\underline{x})}}{|B(\underline{x})|}.$$

Thus, by using the FKG inequality we have that

$$\Phi\left[\widehat{M}_{\widehat{\Pi}_*(M)} \geq \frac{(p - p_c)^\xi}{4} Mn\right] \leq \Phi\left[\frac{1}{|\underline{\widehat{\Pi}}'_*(M)|} \sum_{\underline{x} \in \widehat{\underline{\Pi}}'_*(M)} \frac{\widehat{M}_{B(\underline{x})}}{|B(\underline{x})|} \geq \frac{(p - p_c)^\xi}{28} \Big| E\right],$$

where $E$ is the event that all the dual edges of the boundaries of the blocks $B(\underline{x})$ are *open*. Furthermore, from Proposition 5 of [11], we have that for every $a' \in (\xi + 1, a - \xi)$, there exists a positive constant $c' = c'(a')$ such that for every $\underline{x} \in \widehat{\underline{\Pi}}'_-(M) \cup \widehat{\underline{\Pi}}'_+(M)$, we have

$$m > c'(p - p_c)^{-b'} \quad \Rightarrow \quad \Phi\left[\frac{\widehat{M}_{B(\underline{x})}}{|B(\underline{x})|} \Big| E\right] \leq \Phi^{\widehat{p},w}_{B(\underline{x})}\left[\frac{\widehat{M}_{B(\underline{x})}}{|B(\underline{x})|}\right] \leq \frac{(p - p_c)^\xi}{56}.$$

Observe that the random variables $(\widehat{M}_{B(\underline{x})}/|B(\underline{x})|, \underline{x} \in \widehat{\underline{\Pi}}'_-(M) \cup \widehat{\underline{\Pi}}'_+(M))$ that take their values in $[0, 1]$, are independent and identically distributed under $\Phi[\cdot | E]$. Also, by choosing $c = c'/10$ in (15), we obtain that their mean is bounded above by $(p - p_c)^\xi/56$. Therefore, we can apply Lemma 7 to get

$$\log \Phi\left[\widehat{M}_{\widehat{\Pi}_*(M)} \geq (p - p_c)^\xi \frac{Mn}{4}\right] \leq -\lambda(p - p_c)^{2\xi} \frac{nM}{m^2},$$

where $* = -, +$ and $\lambda > 0$ is a positive constant. $\qquad \square$

5.3.2. *Proof of the interface lemma.* Now we have all the ingredients to give an upper bound on the probability of Sep that captures the existence of an interface.

PROOF OF LEMMA 14.   First we approximate Sep with Sep$'$, and we have

$$(16) \qquad \Phi[\mathrm{Sep}(\delta)] \leq \Phi[\mathrm{Sep}_D(\delta) \setminus \mathrm{Sep}'_D(4\delta)] + \Phi[\mathrm{Sep}'_D(4\delta)].$$

Let $a > 5, \alpha \in (0, (1 + \frac{1}{8a})^{-1})$. Lemma 15 ensures the existence of a positive constant $c_1$ such that for every $M$ satisfying

$$(17) \qquad c_1(p - p_c)^{-a/\alpha} < M < \eta r n/3,$$

we can bound the first term of (16) by

$$(18) \quad \log \Phi[\mathrm{Sep}_D(\delta) \setminus \mathrm{Sep}'_D(4\delta)] \leq -c_1 \frac{n^2}{M^2} \delta(p - p_c)^{1/8} \log\left(\frac{1}{p - p_c}\right).$$

Next, we turn to the estimation of the second term of (16). We fix $\xi > 2$ and decompose the event under consideration as follows:

$$(19) \qquad \begin{aligned} \Phi[\mathrm{Sep}'_D(4\delta)] &\leq \Phi[\exists h\ K(h) \geq (p - p_c)^\xi n] \\ &\quad + \Phi[\{\forall h\ K(h) < (p - p_c)^\xi n\} \cap \mathrm{Sep}'_D(4\delta)], \end{aligned}$$

where $K(h)$ is the number of big open dual clusters that intersect $\widehat{\Pi}(h)$ [$\widehat{\Pi}(h)$ is defined before Lemma 18], and $h$ takes its values in $(-2r\eta n/3, -r\eta n/3) \cup (r\eta n/3, 2r\eta n/3)$. Next, we impose further to the exponent $a$ to be larger than $2\xi + 1$ so that by (17) and by Lemma 18, there exists a positive $\lambda = \lambda(a, \xi)$ such that for $p$ close enough to $p_c$, we have

$$(20) \qquad \Phi[\exists h\ K(h) \geq (p - p_c)^\xi n] \leq n \exp(-\lambda(p - p_c)^{4\xi + 2} nM).$$

Now we turn to the second term of (19). By Lemma 17 we can bound from above $\Phi[\{\forall h\ K(h) < (p - p_c)^\xi n\} \cap \mathrm{Sep}'_D(4\delta)]$ by

$$(21) \qquad \sum_{k=1}^{\lfloor (p-p_c)^\xi n \rfloor} \sum_{\substack{(\widehat{\Lambda}_1, \dots, \widehat{\Lambda}_k) \in \Xi(k) \\ (s_1, \dots, s_k) \in \Upsilon(k)}} \Phi[A(\widehat{\Lambda}_1, \dots, \widehat{\Lambda}_k; s_1, \dots, s_k)],$$

where $\Xi(k)$ and $\Upsilon(k)$ have been defined just before Lemma 17. By Stirling's formula, for all $1 \leq k \leq \lfloor (p - p_c)^\xi n \rfloor$, the cardinality of the sets $\Xi(k)$ and $\Upsilon(k)$ are uniformly bounded from above by

$$(22) \qquad \begin{aligned} |\Xi(k)| &\leq \binom{2n}{k} \leq \exp(4(p - p_c)^\xi n \log n), \\ |\Upsilon(k)| &\leq \binom{2n + k}{k} \leq \exp(4(p - p_c)^\xi n \log n). \end{aligned}$$

Next, we fix $k \in [1, \lfloor (p - p_c)^\xi n \rfloor], (\widehat{\Lambda}_1, \ldots, \widehat{\Lambda}_k) \in \Xi(k)$ and $(s_1, \ldots, s_k) \in \Upsilon(k)$ and we use Lemma 8 to decouple the occurrence of the $k$ separated dual open connections appearing in the event $A(\cdot)$. To do so, we must require that the distance $\ell$ separating the regions $\widehat{\Lambda}_j$'s is large enough. More precisely, there exists a positive constant $c_2$ such that if

$$(23) \qquad \delta n \geq \ell = c_2 (p - p_c)^{-a},$$

then we can apply Lemma 8 $k$ times and use (14) to obtain

$$\Phi[A(\widehat{\Lambda}_1, \ldots, \widehat{\Lambda}_k; s_1, \ldots, s_k)]$$

$$(24) \qquad \leq 2^k \prod_{j=1}^{k} \Phi[\exists \text{ an open dual path } \widehat{\gamma}_j \subset \widehat{\Lambda}_j \text{ of } \operatorname{diam}_w \widehat{\gamma}_j = s_j]$$

$$\leq (2n)^k \exp\left[ -\tau_p(w) n \left( 2\rho - \frac{6\delta\pi}{\eta} r^2 - \frac{2k\ell}{n} \right) \right].$$

Combining (24) and (22) we can bound (21) from above by

$$(25) \qquad \exp\left( -\tau_p(w) n \left( 2\rho - \frac{6\delta\pi}{\eta} r^2 \right) + 8(p - p_c)^\xi n \log n \right)$$

$$\times \sum_{k=1}^{\lfloor n(p-p_c)^\xi \rfloor} \exp(k(\log(2n) + 2\tau_p(w)\ell)).$$

In order to satisfy condition (23), we are limited to regimes where $n \geq c_2 \delta^{-1} (p - p_c)^{-a}$. By making $c_2$ a bit bigger and by using Proposition 4, we can find $n_0 > 1$ such that for every $n > n_0$, we have

$$\Phi[\{ \forall h \ K(h) < (p - p_c)^\xi n \} \cap \operatorname{Sep}'_D(4\delta)]$$

$$(26) \qquad \leq \exp(10(p - p_c)^\xi n \log n + \tau_c c_2 (p - p_c)^{-a+\xi+1} n)$$

$$\times \exp\left( -\tau_c (p - p_c) n \left( 2\rho - \frac{14\delta\pi}{\eta} r^2 \right) \right).$$

By choosing $a > \xi + 1$ and $\log n / \log(1/(p - p_c))$ bounded from above, the first exponential becomes negligible. It remains to specify a regime satisfying (17) and (23) such that the bounds (18) and (20) are smaller than (26). That is, we have to choose $\xi > 2, a > 2\xi + 1, \alpha \in (0, (1 + \frac{1}{8a})^{-1})$ and $n \uparrow \infty, p \downarrow p_c, M$ such that

$$\eta r n / 3 \geq M \geq \max\left( c_1 (p - p_c)^{-a/\alpha}, \right.$$

$$(27) \qquad \left. \left( 2\rho\tau_c (p - p_c) + \frac{\log n}{n} \right) \middle/ (\lambda (p - p_c)^{4\xi+2}) \right),$$

$$n \geq \max\left(\frac{c_2}{\delta}(p-p_c)^{-a}, \frac{2\rho}{\delta c_1}\tau_c M^2 \frac{(p-p_c)^{7/8}}{-\log(p-p_c)}\right).$$

For the choice $\xi > 2, a = 2\xi + 2, \alpha = (1 + \frac{1}{8a})^{-1}, M = (p-p_c)^{4\xi+1}, n = (p - p_c)^{-8\xi-9/8}$, it is easy to check that there exists $n_0 > 1$ such that for all $n > n_0$, the conditions (27) are satisfied. Since $\xi$ has to be larger than 2, we obtain that for every $\gamma > 20$, if $n \uparrow \infty, p \downarrow p_c$ in such a way that $n(p-p_c)^\gamma \to \infty$ then it is possible to find $\xi, a, \alpha, M$ such that asymptotically the conditions (27) are satisfied. Finally, we obtain

$$\limsup_{(n,p)\to(\infty,p_c)} \frac{1}{(p-p_c)n}\log\Phi[\text{Sep}(\delta)] \leq -2\rho\tau_c + \frac{14\delta}{\eta}r^2\tau_c.$$

By choosing $\eta = r\sqrt{2\delta}/14$ and $\rho = r\sqrt{1-\rho}$, we get the desired result.  $\square$

### 5.4. *Potential applications in higher dimensions.*

As promised in the Introduction, we include some comments about the adaptability of our arguments to, say, Bernoulli percolation near criticality in dimensions greater than two. Even if we suppose that an adequate surface tension in a joint limit $p \downarrow p_c$ and $n \uparrow \infty$ exists, the proof of the interface lemma would still need nontrivial adaptation. Indeed, if the critical exponent for the surface tension is greater than the one for the density of the infinite cluster then it is not possible to apply the cutting procedure as in [9]. Furthermore, collecting the surface energy of piecewise interfaces as we did here would be more difficult in higher dimensions. Apart from this obstacle, we believe that the rest of our arguments are robust enough to be adapted in such a situation.

### 5.5. *The local upper bound.*

First we give a result of geometric measure theory. To state this result, let us mention that for every set of finite perimeter $A$, it is possible to define his essential boundary $\partial^* A$ which is a subset of the topological boundary $\partial A$. Also, for each point $x$ of $\partial^* A$, it is possible to define a measure theoretical normal vector $\nu_A(x)$ to $\partial A$ which is taken to point to the exterior of $A$. In order to prove the upper bound we will need the following approximation result.

LEMMA 19. *Let $A$ be a subset of $Q$ having finite perimeter. For any positive $\varepsilon, \delta$, there exists a finite collection of disjoint balls $B(x_i, r_i), i \in I \cup I^{\text{bd}}$, such that:*

- *$\forall i \in I$  $x_i \in \partial^* A \cap \overset{\circ}{Q}$ and $B(x_i, r_i) \subset \overset{\circ}{Q}$;*
- *$\forall i \in I^{\text{bd}}$  $x_i \in \partial^* A \cap \partial Q$ and $B_-(x_i, r_i, \nu_Q(x_i)) \subset Q$;*
- *$\forall i \in I \cup I^{\text{bd}}$  $\mathcal{L}^2((A \cap B(x_i, r_i))\Delta B_-(x_i, r_i, \nu_A(x_i))) \leq \delta r^2$;*

- *Finally*

$$\left| \mathcal{P}(A) - \sum_{i \in I \cup I^{\mathrm{bd}}} 2r_i \right| \le \varepsilon.$$

PROOF.   The proof can be found in Lemma 82 of [10].   □

LEMMA 20.   *Let $\nu \in \mathcal{M}(Q)$ be such that $\mathcal{J}(\nu) < +\infty$. If $\beta \downarrow \beta_c$ and $n \uparrow \infty$ in such a way that $n(\beta - \beta_c)^{20} \to \infty$, and such that $\log(n)/\log(1/(\beta - \beta_c))$ stays bounded, then for every $\varepsilon > 0$ there exists a weak neighborhood $\mathcal{U}$ of $\nu$ in $\mathcal{M}(Q)$ such that*

$$\limsup_{(n,\beta)} \frac{1}{n(\beta - \beta_c)} \log \mu_{\widetilde{\Lambda}(n)}^{+,\beta}(\sigma_n \in \mathcal{U}) \le -(1 - \varepsilon)\mathcal{J}(\nu).$$

PROOF.   By definition of $\mathcal{J}$, since $\mathcal{J}(\nu) < +\infty$, there exists a Borel subset $A$ of $Q$ such that $\nu$ is the measure with density $-1_A + 1_{Q \setminus A}$ with respect to the Lebesgue measure, and

$$\mathcal{J}(\nu) = \tau_c \mathcal{P}(A) = \tau_c \mathcal{P}(Q \setminus A).$$

If $\mathcal{P}(A) = 0$, there is nothing to prove. Suppose that $\mathcal{P}(A) > 0$. Let $\delta_0, \varepsilon' \in \,]0,1[\,$ that we will tune later. Let $B(x_i, r_i), i \in I \cup I^{\mathrm{bd}}$, be a finite collection of disjoint balls associated with $A, \varepsilon'$ and $\delta_0/3$, as given in Lemma 19. For $i$ in $I \cup I^{\mathrm{bd}}$, let $f_i, g_i$ be two continuous functions having compact support and taking values in $[0, 1]$ such that

$$\forall x \in \mathbb{R}^2 \setminus \mathring{B}_-(x_i, r_i, \nu_A(x_i)) \qquad f_i(x) = 0,$$

$$\forall x \in \overline{B}_+(x_i, r_i, \nu_A(x_i)) \qquad g_i(x) = 1,$$

$$\left(\frac{1}{2} - \frac{\delta_0}{4}\right)\pi r_i^2 \le \int f_i \, d\mathcal{L}^2, \qquad \int g_i \, d\mathcal{L}^2 \le \left(\frac{1}{2} + \frac{\delta_0}{4}\right)\pi r_i^2,$$

where $\nu_A(x_i)$ is the exterior normal vector of $A$ at $x_i$ and where $\mathring{B}_-$ and $\overline{B}_+$ denote the interior and the closure of the half balls. Also, we require that there exists $s_i > 0$ such that, if we set

$$D_-^i = \{y \in B_-(x_i, r_i, \nu_A(x_i)) : \mathrm{d}(y, \mathbb{R}^2 \setminus B_-(x_i, r_i, \nu_A(x_i))) \le s_i\},$$

$$D_+^i = \{y \notin B_+(x_i, r_i, \nu_A(x_i)) : \mathrm{d}(y, B_+(x_i, r_i, \nu_A(x_i))) \le s_i\},$$

then we have

$$\mathcal{L}^2(D_-^i) \le \frac{\delta_0}{8}\pi r_i^2, \qquad \forall x \in B_-(x_i, r_i, \nu_A(x_i)) \setminus D_-^i \qquad f_i(x) = 1,$$

$$\mathcal{L}^2(D_+^i) \le \frac{\delta_0}{8}\pi r_i^2, \qquad \forall x \in \mathbb{R}^d \setminus B_+(x_i, r_i, \nu_A(x_i)) \setminus D_+^i \qquad g_i(x) = 0.$$

These conditions imply that

$$\nu(f_i) = -\int_A f_i \, d\mathcal{L}^2 + \int_{Q \setminus A} f_i \, d\mathcal{L}^2$$

$$\leq -\int f_i \, d\mathcal{L}^2 + 2\mathcal{L}^2(B_-(x_i, r_i, \nu_A(x_i)) \setminus A)$$

$$\leq (-1 + \delta_0)\frac{1}{2}\pi r_i^2,$$

$$\nu(g_i) = -\int_A g_i \, d\mathcal{L}^2 + \int_{Q \setminus A} g_i \, d\mathcal{L}^2$$

$$\geq \int g_i \, d\mathcal{L}^2 - 2\mathcal{L}^2(A \cap B_+(x_i, r_i, \nu_A(x_i))) - 2\mathcal{L}^2(D_+^i).$$

$$\geq (1 - \delta_0)\frac{1}{2}\pi r_i^2.$$

Let $\mathcal{U}$ be the weak neighborhood of $\nu$ in $\mathcal{M}(Q)$ defined by

$$\mathcal{U} = \left\{ \varrho \in \mathcal{M}(Q) : \forall i \in I \; \rho(f_i) < \nu(f_i) + \frac{\delta_0}{2}\pi r_i^2, \rho(g_i) > \nu(g_i) - \frac{\delta_0}{2}\pi r_i^2 \right\}.$$

Next, we choose $\gamma > 20$ and consider $n \uparrow \infty, p \downarrow p_c$ in such a way that $n(p - p_c)^\gamma \to \infty$. We rescale the lattice by a factor $M = \sqrt{n}$ and choose $L = \sqrt{M}$. Let $\delta > 0$, for $\underline{x} \in \mathbb{Z}^2$, the block variable $X(\underline{x})$ is defined as the indicator function of the event

$$R(B'(\underline{x}), L) \cap V(B(\underline{x}), \delta) \cap W(B(\underline{x}), \delta) \cap T(B(\underline{x}), \delta).$$

Let us fix $i \in I$. Let $*$ be a symbol representing either $-$ or $+$. We define

$$B_*(n, i) = B_*(nx_i, nr_i, \nu_A(x_i)),$$

$$\underline{B}_*(n, i) = \{\underline{x} \in \mathbb{Z}^2 : B(\underline{x}) \subset \mathring{B}_*(n, i)\}.$$

By the above choice of $n, p$ and $M$, for $n$ large enough we have

$$\mathcal{L}^2\left(B_*(n, i) \setminus \bigcup_{\underline{x} \in \underline{B}_*(n,i)} B(\underline{x})\right) \leq \theta \delta \mathcal{L}^2(B_*(n, i)).$$

We define $\mathcal{S}_*$ as the collection of the clusters which are included in one of the boxes $B(\underline{x}), \underline{x} \in \underline{B}_*(n, i)$, but which do not intersect their boundaries;

$$\mathcal{S}_* = \bigcup_{\underline{x} \in \underline{B}_*(n,i)} \{C \text{ cluster in } B(\underline{x}) \text{ such that } C \cap \partial B(\underline{x}) = \varnothing\}.$$

Let $\underline{x} \in \underline{B}_*(n, i)$ be such that $X(\underline{x}) = 1$. Then

$$\sum_{C \in \mathcal{S}_*} |C \cap B(\underline{x})| \geq |B(\underline{x})| - |\{x \in B(\underline{x}) : x \leftrightarrow \partial B(\underline{x})\}| \geq (1 - \theta - \delta\theta)|B(\underline{x})|.$$

Summing these inequalities, we get

$$\sum_{C \in \mathcal{S}_*} |C \cap B_*(n,i)| \geq (1 - \theta - 2\delta\theta)\mathcal{L}^2(B_*(n,i)). \tag{28}$$

We define also $\mathcal{C}$ as the collection of the $B(n,i)$ clusters which do not belong to $\mathcal{S}_+ \cup \mathcal{S}_-$. For a cluster $C$, we denote by $\sigma(C)$ its color. For $\underline{x} \in \underline{B}_-(n,i) \cup \underline{B}_+(n,i)$, we have

$$\theta n^2 \sigma_n(B(\underline{x})) = \sum_{C \in \mathcal{C}} \sigma(C)|C \cap B(\underline{x})| + \sum_{\substack{x \in B(\underline{x}) \\ x \not\leftrightarrow \partial B}} \sigma(x).$$

Whenever $X(\underline{x}) = 1$, the event $T$ occurs and the modulus of the last sum is less than $\delta\theta|B(\underline{x})|$. Suppose that, for $* = -$ and $* = +$

$$M^2 \sum_{\underline{x} \in \underline{B}_*(n,i)} 1_{X(\underline{x})=0} \leq \delta\theta\mathcal{L}^2(B_*(n,i)).$$

Summing the previous inequalities, we get

$$\theta n^2 \sigma_n(\mathring{\bar{B}}_+(n,i)) \leq \sum_{C \in \mathcal{C}} \sigma(C)|C \cap B_+(n,i)| + 5\delta\theta\mathcal{L}^2(B_+(n,i)),$$

$$\theta n^2 \sigma_n(\overline{B}_-(n,i)) \geq \sum_{C \in \mathcal{C}} \sigma(C)|C \cap B_-(n,i)| - 5\delta\theta\mathcal{L}^2(B_-(n,i)).$$

Let us denote by $\mathcal{C}_-$ (respectively, $\mathcal{C}_+$) the collection of the negatively (respectively, positively) colored clusters of $\mathcal{C}$. Notice that the collections $\mathcal{S}_+ \cup \mathcal{S}_-, \mathcal{C}_-, \mathcal{C}_+$ are disjoint. Suppose in addition that $\sigma_n \in \mathcal{U}$. The very definition of the neighborhood $\mathcal{U}$, the two previous inequalities and (28) yield that

$$\sum_{C \in \mathcal{C}_- \cup \mathcal{S}_-} |C \cap B_+(n,i)| \leq \theta(8\delta + \delta_0)\mathcal{L}^2(B_+(n,i)),$$

$$\sum_{C \in \mathcal{C}_- \cup \mathcal{S}_-} |C \cap B_-(n,i)| \geq (1 - \theta(8\delta + \delta_0))\mathcal{L}^2(B_-(n,i)).$$

Thus the collection $\mathcal{C}_- \cup \mathcal{S}_-$ realizes the event $\text{Sep}(n, x_i, r_i, w_i, 8\delta + \delta_0)$. In fact, some care is needed on the boundary of $Q$, and one needs to define a variant of the event Sep for the balls that intersect the boundary. As the same reasoning holds in this situation, we omit the details. Choosing $\delta < \delta_0/8$, we conclude that

$$\mu_{\Lambda(n)}^{+,\beta}[\sigma_n \in \mathcal{U}] \leq \sum_{*=-,+} \sum_{i \in I \cup I^{\text{bd}}} \mathbb{P}_n^+\left(M^2 \sum_{\underline{x} \in \underline{B}_*(n,i)} 1_{X(\underline{x})=0} > \theta\delta\mathcal{L}^2(B_*(n,i))\right)$$

$$+ \Phi_{\Lambda(n)}^{w,p}\left[\bigcap_{i \in I \cup I^{\text{bd}}} \text{Sep}(n, x_i, r_i, \nu_A(x_i), 2\delta_0)\right].$$

Since we chose $\gamma > 20$, $p \downarrow p_c$, $n \uparrow \infty$ in such a way that $n(p - p_c)^\gamma \to \infty$, $M = \sqrt{n}, L = \sqrt{M}$, we can apply Lemmas 9, 10, 11, 12 to conclude that the block process $(X(\underline{x}))$ satisfies the hypothesis of Lemma 6 and that for all $i \in I \cup I^{\mathrm{bd}}$ and for $* = +, -$ the mean of the normalized sum below converges faster to zero than $\theta$. Thus

$$\limsup_{(n,p)} \frac{1}{n(p - p_c)} \log \mathbb{P}_n^+\left[\frac{M^2}{\mathcal{L}^2(B_*(n,i))} \sum_{\underline{x} \in \underline{B}_*(n,i)} 1_{X(\underline{x})=0} > \delta\theta\right] = -\infty.$$

Next, the sets $B(x_i, r_i), i \in I \cup I^{\mathrm{bd}}$, are compact and disjoint. Also $I \cup I^{\mathrm{bd}}$ is finite and fixed, thus applying Lemma 14 with Lemma 8, we get

$$\limsup_{(n,\beta)} \frac{1}{(\beta - \beta_c)n} \log \mu_{\Lambda(n)}^{+,\beta}[\sigma_n \in \mathcal{U}] \leq -\sum_{i \in I \cup I^{\mathrm{bd}}} 2 r_i \tau_c (1 - c'\sqrt{2\delta_0})$$

$$\leq -\tau_c(\mathcal{P}(A) - \varepsilon')(1 - c'\sqrt{2\delta_0}),$$

where $c'$ is the constant appearing in Lemma 14. Let $\varepsilon > 0$. By choosing $\varepsilon'$ such that $\varepsilon = \varepsilon'(1 + 1/\mathcal{P}(A))$ and $\delta_0$ such that $c'\sqrt{2\delta_0} < \varepsilon'$, we finally get the desired upper bound.  □

### 5.6. Exponential tightness.

In order to prove the exponential tightness, we proceed as in [10]. The same approach works in our context with the exception that one has to be careful with the scales of renormalization and some extra care is needed because $\theta$ converges to 0 when $p \downarrow p_c$. As in [10], we will first define a roughening $\widetilde{\sigma}_n$ of the random measure $\sigma_n$. This auxiliary measure will be regular enough to produce suitable surface energy estimates, and the proof will be completed by proving that the two random measures $\sigma_n$ and $\widetilde{\sigma}_n$ are exponentially contiguous.

### 5.6.1. The rough measure $\widetilde{\sigma}_n$ and surface energy estimates.

In order to construct $\widetilde{\sigma}_n$ we will work with the box $\Lambda(n)$ rescaled by a factor $K$ that will depend on $p$ in a way to be made precise in the course of our analysis. The renormalized box will be denoted by

$$\underline{\Lambda}(n) = \{\underline{x} \in \mathbb{Z}^2 : B(\underline{x}) \cap \Lambda(n) \neq \varnothing\}.$$

On $\underline{\Lambda}(n)$ we define the 0–1 valued random field $X(\underline{x})$, $\underline{x} \in \underline{\Lambda}(n)$, by

$$\forall \underline{x} \in \underline{\Lambda}(n) \qquad X(\underline{x}) = 1_{R(B'(\underline{x}),K)},$$

where $R(B'(\underline{x}),K)$ has been defined just before Lemma 9. When $X(\underline{x}) = 1$ we will say that the block $B(\underline{x})$ is good, and if $X(\underline{x}) = 0$ the block will be said to be bad. We also need a very similar filling procedure as before but in the renormalized lattice. Namely, for every connected subset $\underline{A}$ of $\underline{\Lambda}(n)$ we define

$$\text{fill}\,\underline{A} = \underline{A} \cup \bigcup \underline{R},$$

where the union runs over the residual $\mathbb{L}^\infty$-connected components $\underline{R}$ of $\underline{A}$ such that $\operatorname{diam}_\infty \underline{R} < \log n$ and $\underline{R} \cap \partial\Lambda(n) = \varnothing$. For each cluster $C$, we define its coarse graining $\underline{C} = \{\underline{x} \in \underline{\Lambda}(n) : B(\underline{x}) \cap C \neq \varnothing\}$. We say that a cluster $C$ is large if $\operatorname{diam} C \geq K \log n$, in which case we have $\operatorname{diam} \underline{C} \geq \log n$.

For a large cluster $C$ we define

$$\widehat{\underline{C}} = \bigcup \text{ fill } \underline{A},$$

where the union runs over all the connected components $\underline{A}$ of good blocks such that $\underline{C} \cap \underline{A} \neq \varnothing$. Lemma 18.1 of [10] ensures that if $C_1, C_2$ are two distinct large clusters, then $\widehat{\underline{C}}_1 \cap \widehat{\underline{C}}_2 = \varnothing$. We define $\widetilde{\sigma}_n$ as the random measure on $Q$ whose density with respect to the Lebesgue measure is $1_{P_n} - 1_{M_n} + 1_{(Q\setminus P_n)\setminus M_n}$ where

$$P_n = \bigcup_{\substack{C \text{ large cluster} \\ \sigma(C)=+}} \bigcup_{\underline{x}\in\widehat{\underline{C}}} \frac{1}{n}B(\underline{x}) \cup \left( Q \setminus \Lambda\left(1 - \frac{6K}{n}\right) \right),$$

$$M_n = \bigcup_{\substack{C \text{ large cluster} \\ \sigma(C)=-}} \bigcup_{\underline{x}\in\widehat{\underline{C}}} \frac{1}{n}B(\underline{x}) \setminus \left( Q \setminus \Lambda\left(1 - \frac{6K}{n}\right) \right).$$

This measure $\widetilde{\sigma}_n$ is regular enough to establish the required exponential tightness. To do this, we consider the set $\widehat{\underline{F}}$ of all the $\mathbb{L}^\infty$-connected components of bad blocks in $\underline{\Lambda}(n)$ that intersect $\underline{\Lambda}(n - 6K \log n)$ and one of the sets,

$$\partial_\infty^{\text{out}}\widehat{\underline{C}} = \{\underline{x} \notin \widehat{\underline{C}} : \exists \underline{y} \in \widehat{\underline{C}}, |\underline{x} - \underline{y}|_\infty = 1\},$$

where $C$ is a large cluster.

LEMMA 21.   *Let $a > 5$. There exist two positive constants $\lambda$ and $c = c(a)$ such that if $K > c(p - p_c)^{-a}$ and $n$ is such that $n > K \log n$, then, for $u > 0$,*

$$\Phi_{\Lambda(n)}^w[|\widehat{\underline{F}}| \geq u] \leq \exp(-\lambda(p - p_c)Ku)$$

*and*

$$\mathbb{P}_n^+[\mathcal{P}(P_n) + \mathcal{P}(M_n) \geq u] \leq \exp\left( -\lambda\frac{u - 16}{8}(p - p_c)n \right).$$

PROOF.   Let $C$ be a large cluster. By definition, the $\mathbb{L}^\infty$ outer boundary of $\widehat{\underline{C}}$ consists of bad blocks whenever $\widehat{\underline{C}} \neq \varnothing$. In the case $\widehat{\underline{C}} = \varnothing$, we define $\partial_\infty^{\text{out}}\widehat{\underline{C}}$ as $\widehat{\underline{C}}$ which again consists only of bad blocks. Let $\underline{F}$ be an $\mathbb{L}^\infty$-connected component of bad blocks intersecting $\partial_\infty^{\text{out}}\widehat{\underline{C}}$ and the rescaled box $\underline{\Lambda}(n - 6K \log n)$. As proved in Lemma 18.2 of [10], we have $|\underline{F}| \geq \log n$. Thus

$$\Phi_{\Lambda(n)}^w[|\widehat{\underline{F}}| \geq u] \leq \Phi_{\Lambda(n)}^w[|\{\underline{x} \in \underline{\Lambda}(n - 6K \log n) : |\underline{C}(\underline{x})| \geq \log n\}| \geq u],$$

where $\underline{C}(\underline{x})$ denotes the $\mathbb{L}^\infty$-connected component of occupied sites containing $\underline{x}$ [a site $\underline{y}$ is said to be occupied if $X(\underline{y}) = 0$]. Let $a > 5$, Lemma 9 ensures the existence of two positive constants $\lambda$ and $c = c(a)$ such that if $K > c(p - p_c)^{-a}$, then

$$\sup_{\underline{x} \in \underline{\Lambda}(n)} \log \Phi[X(\underline{x}) = 0 | X(\underline{z}), |\underline{x} - \underline{z}|_\infty \geq 3] \leq -\lambda(p - p_c)K.$$

Combining this with Lemma 5, we get that

$$\Phi^w_{\Lambda(n)}[|\widehat{\underline{F}}| \geq u] \leq 2 \sum_{j \geq u} \exp\left(j\left(c - \frac{1}{9}\lambda(p - p_c)K\right)\right),$$

which yields the first claim.

For the second claim, note that the boundaries of $P_n$ and $M_n$ are located either on $\partial Q \cup \partial \Lambda(1 - 6K \log n/n)$ or on the faces of the blocks of $\widehat{\underline{F}}$. Thus $\mathcal{P}(P_n) + \mathcal{P}(M_n) \leq 16 + 8\frac{K}{n}|\widehat{\underline{F}}|$ so that

$$\mathbb{P}^+_n[\mathcal{P}(P_n) + \mathcal{P}(M_n) \geq u] \leq \Phi^w_{\Lambda(n)}\left[|\widehat{\underline{F}}| \geq (u - 16)\frac{n}{8K}\right]$$

$$\leq \exp\left(-(p - p_c)\lambda K(u - 16)\frac{n}{8K}\right),$$

which yields the desired result. $\quad\square$

5.6.2. *Exponential contiguity.* Here we show that the rough measure $\widetilde{\sigma}_n$ is a good approximation of the original random measure $\sigma_n$. Let $f : \mathbb{R}^2 \to \mathbb{R}$ be a continuous function having compact support. In order to estimate $|\sigma_n(f) - \widetilde{\sigma}_n(f)|$, we use another block coarse graining with scale $L = K \log n$. We fix $\varepsilon > 0$ and define for $y \in \mathbb{Z}^2$ the block variable $Y(\underline{y})$ as the indicator function of the event

$$R(B'(\underline{y}), \sqrt{L}) \cap V(B(\underline{y}), \varepsilon) \cap W(B(\underline{y}), \varepsilon).$$

Let

$$\underline{A} = \{\underline{y} \in \mathbb{Z}^2 : B(\underline{y}) \cap \Lambda(n - 6L) \neq \varnothing\}.$$

Note that $|\underline{A}|L^2 \leq n^2$. We further introduce

$$\forall \underline{x} \in \underline{\Lambda}(n) \qquad B_n(\underline{x}) = \frac{1}{n}B(\underline{x}),$$

$$A_0 = \bigcup_{\underline{y} \in \underline{A} : Y(\underline{y}) = 0} B_n(\underline{y}),$$

$$A_1 = \bigcup_{\underline{y} \in \underline{A} : Y(\underline{y}) = 1} B_n(\underline{y})$$

and make the following decomposition:

$$|\sigma_n(f) - \widetilde{\sigma}_n(f)| \leq \left| \int_{A_0} f(x)(d\sigma_n(x) - d\widetilde{\sigma}_n(x)) \right|$$

(29)
$$+ \left| \int_{(Q \setminus A_0) \setminus A_1} f(x)(d\sigma_n(x) - d\widetilde{\sigma}_n(x)) \right|$$

$$+ \left| \int_{A_1} f(x)(d\sigma_n(x) - d\widetilde{\sigma}_n(x)) \right|.$$

We bound the first term of (29) as follows:

$$\left| \int_{A_0} f(x)(d\sigma_n(x) - d\widetilde{\sigma}_n(x)) \right| \leq 2\|f\|_\infty \frac{1}{\theta |\underline{A}|} \sum_{\underline{y} \in \underline{A}} 1_{Y(\underline{y})=0}.$$

In a similar way we bound the second term:

$$\left| \int_{(Q \setminus A_0) \setminus A_1} f(x)(d\sigma_n(x) - d\widetilde{\sigma}_n(x)) \right|$$

$$\leq \frac{\mathcal{L}^2(\Lambda(n) \setminus \Lambda(n-6L))}{n^2} \left( \frac{1}{\theta} + 1 \right) \|f\|_\infty \leq \frac{24L}{\theta n} \|f\|_\infty.$$

For the third term of (29), we further decompose it into

$$\left| \int_{A_1} f(x)(d\sigma_n(x) - d\widetilde{\sigma}_n(x)) \right|$$

$$\leq \frac{1}{\theta n^2} \left| \sum_{\underline{y} \in \underline{A} \,:\, Y(\underline{y})=1} \sum_{\substack{C \subset B(\underline{y}) \\ C \cap \partial^{\mathrm{in}} B(\underline{y}) = \varnothing}} \sigma(C) \sum_{x \in C} f(x/n) \right|$$

(30)
$$+ \frac{1}{\theta n^2} \left| \sum_{\underline{y} \in \underline{A} \,:\, Y(\underline{y})=1} \sum_{\substack{x \in B(\underline{y}) \setminus C(\underline{y}) \\ x \leftrightarrow \partial^{\mathrm{in}} B(\underline{y})}} \sigma(x) f(x/n) \right|$$

$$+ \left| \frac{1}{\theta n^2} \sum_{\substack{\underline{y} \in \underline{A} \\ Y(\underline{y})=1}} \sigma(C(\underline{y})) \sum_{x \in C(\underline{y})} f(x/n) - \int_{A_1} f(x) \, d\widetilde{\sigma}_n(x) \right|,$$

where $C(\underline{y})$ denotes the unique crossing cluster of $B(\underline{y})$ whenever $Y(\underline{y})=1$. Using the definition of the good blocks, we further bound the second term of the right-hand side in (30):

$$\frac{1}{\theta n^2} \left| \sum_{\underline{y} \in \underline{A} \,:\, Y(\underline{y})=1} \sum_{\substack{x \in B(\underline{y}) \setminus C(\underline{y}) \\ x \leftrightarrow \partial^{\mathrm{in}} B(\underline{y})}} \sigma(x) f(x/n) \right| \leq 2\varepsilon \|f\|_\infty.$$

Since $f$ is continuous and has compact support, we have for $K \log n/n$ small enough that for all $y \in \mathbb{Z}^2 \sup_{x,z \in B_n(\underline{y})} |f(x) - f(z)| \le \varepsilon \|f\|_\infty$. Using this observation and the properties of the good blocks, we can bound the third term of (30) by

$$
\left| \frac{1}{\theta n^2} \sum_{\substack{\underline{y} \in \underline{A} \\ Y(\underline{y})=1}} \sigma(C(\underline{y})) \sum_{x \in C(\underline{y})} f(x/n) - \int_{A_1} f(x)\, d\widetilde{\sigma}_n(x) \right|
$$

$$
\le \left| \left[ \sum_{\substack{\underline{y} \in \underline{A} \\ Y(\underline{y})=1}} \left( \frac{1}{\theta n^2} \sigma(C(\underline{y})) |C(\underline{y})| - \widetilde{\sigma}_n(B_n(\underline{y})) \right) \right] \max_{x \in B_n(\underline{y})} f(x) \right|
$$

$$
(31) \qquad + \frac{1}{\theta n^2} \left| \sum_{\substack{\underline{y} \in \underline{A} \\ Y(\underline{y})=1}} \sigma(C(\underline{y})) \sum_{x \in C(\underline{y})} \left( f\left(\frac{x}{n}\right) - \max_{z \in B_n(\underline{y})} f(z) \right) \right|
$$

$$
+ \left| \sum_{\substack{\underline{y} \in \underline{A} \\ Y(\underline{y})=1}} \int_{B_n(\underline{y})} f(x) - \max_{z \in B_n(\underline{y})} f(z)\, d\widetilde{\sigma}_n(x) \right|
$$

$$
\le \|f\|_\infty \sum_{\substack{\underline{y} \in \underline{A} \\ Y(\underline{y})=1}} \left| \frac{1}{\theta n^2} \sigma(C(\underline{y})) |C(\underline{y})| - \widetilde{\sigma}_n(B_n(\underline{y})) \right| + \varepsilon(\varepsilon + 2) \|f\|_\infty.
$$

Next, we study the sum in the above inequality. Let $\underline{y} \in \underline{A}$ such that $Y(\underline{y}) = 1$. Several cases arise.

- $B_n(\underline{y}) \subset M_n$ and $\sigma(C(\underline{y})) = -1$. We have then

$$
\left| \frac{1}{\theta n^2} \sigma(C(\underline{y})) |C(\underline{y})| - \widetilde{\sigma}_n(B_n(\underline{y})) \right| = \frac{1}{n^2} \left| -\frac{|C(\underline{y})|}{\theta} + |B(\underline{y})| \right| \le \varepsilon \frac{|B(\underline{y})|}{n^2}.
$$

- $B_n(\underline{y}) \subset P_n$ and $\sigma(C(\underline{y})) = 1$: this case is symmetric to the previous one and

$$
\left| \frac{1}{\theta n^2} \sigma(C(\underline{y})) |C(\underline{y})| - \widetilde{\sigma}_n(B_n(\underline{y})) \right| \le \varepsilon \frac{|B(\underline{y})|}{n^2}.
$$

- $B_n(\underline{y}) \cap M_n \ne \varnothing$ and $B_n(\underline{y}) \not\subset M_n$, then $B_n(\underline{y})$ meets the boundary of $M_n$.
- $B_n(\underline{y}) \cap P_n \ne \varnothing$ and $B_n(\underline{y}) \not\subset P_n$, then $B_n(\underline{y})$ meets the boundary of $P_n$.

In the two last cases, we bound crudely as follows:

$$
\left| \frac{1}{\theta n^2} \sigma(C(\underline{y})) |C(\underline{y})| - \widetilde{\sigma}_n(B_n(\underline{y})) \right| \le (2 + \varepsilon) \frac{|B(\underline{y})|}{n^2}.
$$

This will suffice because in this case
$$B_n(\underline{y}) \subset \left\{ x \in \mathbb{R}^2 : d_\infty(x, (\partial P_n \cup \partial M_n) \cap \Lambda(1 - 3L/n)) \le \frac{L}{n} \right\}.$$

Moreover,
$$\mathcal{L}^2\left(\left\{ x \in \mathbb{R}^2 : d_\infty(x, (\partial P_n \cup \partial M_n) \cap \Lambda(1 - 3L/n)) \le \frac{L}{n} \right\}\right)$$

$$\le (|\partial^{\mathrm{in}} n P_n| + |\partial^{\mathrm{in}} n M_n|)\left(\frac{2L+2}{n}\right)^2 \le (\mathcal{P}(M_n) + \mathcal{P}(P_n))\frac{9L^2}{n}.$$

- $(B_n(\underline{y}) \subset M_n$ and $\sigma_n(B_n(\underline{y})) > 0)$ or $(B_n(\underline{y}) \subset P_n$ and $\sigma_n(B_n(\underline{y})) < 0)$. These conditions imply that the whole block $B_n(\underline{y})$ has been added to $M_n$ or $P_n$ by the filling operation. Yet the regions which are added by the filling operation have a diameter at most $K(\log n - 1)$, so this case cannot occur.

- $B_n(\underline{y}) \cap M_n = \varnothing$ and $B_n(\underline{y}) \cap P_n = \varnothing$. In this case, $B(\underline{y}) \subset \widehat{F}$.

Summing the previous inequalities, we get

$|\sigma_n(f) - \widetilde{\sigma}_n(f)|$

$$\le \frac{1}{\theta n^2}\left|\sum_{C \in \mathcal{S}} \sigma(C) \sum_{x \in C} f(x/n)\right|$$

$$+ \|f\|_\infty \left(\frac{2}{\theta|\underline{A}|}\sum_{\underline{y} \in \underline{A}} 1_{Y(\underline{y}) = 0} + \frac{L^2}{\theta n^2}|\widehat{F}| + 9(2 + \varepsilon)\frac{L^2}{n}(\mathcal{P}(M_n) + \mathcal{P}(P_n))\right)$$

$$+ \|f\|_\infty \left(\frac{24L}{\theta n} + \varepsilon(\varepsilon + 6)\right),$$

where
$$\mathcal{S} = \bigcup_{\underline{y} \in \underline{A} : Y(\underline{y}) = 1} \{C \in B(\underline{y}) : C \cap \partial^{\mathrm{in}} B(\underline{y}) = \varnothing\}.$$

Note that by the definition of a good block, any cluster $C$ of a good block $B(\underline{y})$ that has a diameter larger than $\sqrt{L} = \sqrt{K \log n}$ is connected to the crossing cluster of $B'(\underline{y})$, and thus, such a cluster $C$ is connected to $\partial^{\mathrm{in}} B(\underline{y})$. Therefore, any cluster of $\mathcal{S}$ has a diameter that is smaller than $\sqrt{K \log n}$. Next, we analyze the deviations of the first term in the last inequality.

Since $|\mathcal{S}| \le n^2$, we have

$$\mathbb{P}_n^+\left[\frac{1}{\theta n^2}\left|\sum_{C \in \mathcal{S}} \sigma(C) \sum_{x \in C} f\left(\frac{x}{n}\right)\right| > \varepsilon \|f\|_\infty\right]$$

$$\le \sum_{\omega \in \Omega_{\Lambda(n)}} \mathbb{P}_n^+\left[\frac{1}{|\mathcal{S}|}\left|\sum_{C \in \mathcal{S}} Y_C\right| > \frac{\varepsilon\theta}{K \log n}\Big|\omega\right]\Phi_{\Lambda(n)}^{w,p}[\omega],$$

where $Y_C = \sigma(C) \sum_{x \in C} f(x/n) / (\|f\|_\infty K \log n)$.

Under $\mathbb{P}_n^+[\cdot|\omega]$, the sequence $(Y_C, C \in \mathcal{S}(\omega))$ is independent and takes its values in $[-1, 1]$ (recall that the diameters of the clusters of $\mathcal{S}$ are bounded by $\sqrt{K \log n}$). Therefore we can apply Theorem 7 to get

$$\mathbb{P}_n^+\left[\frac{1}{\theta n^2}\left|\sum_{C \in \mathcal{S}} \sigma(C) \sum_{x \in C} f\left(\frac{x}{n}\right)\right| > \varepsilon \|f\|_\infty\right]$$

$$\leq 2 \sum_{\omega \in \Omega_{\Lambda(n)}} \exp\left[-\left(\frac{\varepsilon\theta}{K \log n}\right)^2 |\mathcal{S}|\right] \Phi_{\Lambda(n)}^{w,p}[\omega]$$

$$\leq 2\Phi_{\Lambda(n)}^{w,p}\left[\frac{1}{|\underline{A}|}\sum_{\underline{y} \in \underline{A}} 1_{Y(\underline{y})=0} \geq \varepsilon\right]$$

$$+ 2\exp\left(-\frac{(1-\varepsilon)(1-(1+\varepsilon)\theta)\varepsilon^2\theta^2 n^2}{(K \log n)^3}\right),$$

where we used the fact that

$$|\mathcal{S}| \geq \frac{|\{\underline{y} \in \underline{A} : Y(\underline{y}) = 1\}|L^2(1-(1+\varepsilon)\theta)}{\sup_{C \in \mathcal{S}}(\mathrm{diam}^2(C))}.$$

Therefore, by imposing

$$(32) \qquad n \geq \frac{24K}{\varepsilon\theta} \log n,$$

we have

$$\mu_n^{+,\beta}[|\sigma_n(f) - \widetilde{\sigma}_n(f)| > \|f\|_\infty(\varepsilon + 10\varepsilon)]$$

$$(33) \qquad \begin{aligned} &\leq \Phi_{\Lambda(n)}^{w,p}\left[|\widehat{F}| > \varepsilon\theta\frac{n^2}{L^2}\right] + \mathbb{P}_n^+\left[\mathcal{P}(M_n) + \mathcal{P}(P_n) > \varepsilon\frac{n}{L^2}\right] \\ &\quad + 2\exp\left(-\frac{(1-\varepsilon)(1-(1+\varepsilon)\theta)\varepsilon^2\theta^2 n^2}{(K \log n)^3}\right) \\ &\quad + 3\Phi_{\Lambda(n)}^{w,p}\left[\frac{1}{|\underline{A}|}\sum_{\underline{y} \in \underline{A}} 1_{X(\underline{y})=0} \geq \frac{\varepsilon\theta}{2}\right]. \end{aligned}$$

By Lemma 21, for any $a > 5$, there exists a positive constant $c$ such that if we impose that

$$(34) \qquad K > c(p-p_c)^{-a} \quad \text{and} \quad n > K \log n,$$

then

$$(35) \qquad \Phi_{\Lambda(n)}^{w,p}\left[|\widehat{F}| > \varepsilon\theta\frac{n^2}{L^2}\right] \leq \exp\left(-\lambda(p-p_c)\frac{\varepsilon\theta}{K}\frac{n^2}{\log^2 n}\right).$$

If we further impose that

$$(36) \qquad \frac{\varepsilon n}{(K \log n)^2} > 32,$$

then Lemma 21 gives

$$(37) \qquad \mathbb{P}_n^+ \left[ (\mathcal{P}(M_n) + \mathcal{P}(P_n)) > \varepsilon \frac{n}{L^2} \right] \leq \exp \left( -\frac{\lambda \varepsilon}{16} \frac{p - p_c}{K^2} \frac{n^2}{\log^2 n} \right).$$

Furthermore, under condition (34), Lemmas 9, 10, 11 imply that the block process $(Y(\underline{y}), \underline{y} \in \underline{A})$ satisfies the hypothesis of Lemma 6 and that the mean of the normalized sum below converges to zero faster than $\theta$. Thus

$$(38) \qquad \limsup_{n \uparrow \infty, p \downarrow p_c} \frac{1}{(p - p_c)} \log \Phi_{\Lambda(n)}^{w,p} \left[ \frac{1}{|\underline{A}|} \sum_{\underline{y} \in \underline{A}} 1_{Y(\underline{y})=0} \geq \frac{\varepsilon \theta}{2} \right] = -\infty.$$

Finally, we verify that if we choose $a > 5$ and $c$ large enough and impose

$$(39) \qquad K > c(p - p_c)^{-a} \quad \text{and} \quad n > K^3 \log^3 n,$$

then for $p$ close enough to $p_c$, conditions (32), (34) and (36) are satisfied. Moreover, if $n \uparrow \infty$ and $p \downarrow p_c$ in a regime where (39) is satisfied, then (35), (37) and (38) imply that all the terms of (33) decay exponentially fast with a speed larger than $(p - p_c)n$. Thus, when $n \uparrow \infty$ and $p \downarrow p_c$ in such a way that $n(p - p_c)^{20} \uparrow \infty$, we have for every positive $\varepsilon$

$$\limsup_{(n,p)} \frac{1}{n(p - p_c)} \log \mathbb{P}_n^+ (|\sigma_n(f) - \widetilde{\sigma}_n(f)| > \varepsilon) = -\infty.$$

Together with Lemma 21, this concludes the proof of the exponential tightness.

## 6. The lower bound.
We are only left with the lower bound in order to finish our large deviation principle.

### 6.1. *Preparatory lemmas.*
To prove the lower bound we will consider an event whose probability is of the correct order. For this we will use the following lemma to approximate sets of finite perimeter with polyhedral sets.

**LEMMA 22.** *Let $A$ be a subset of $Q = [-1, 1]^2$ having finite perimeter. For any $\varepsilon > 0$, there exists a finite union of polyhedral sets $D$ such that*

$$\overline{D} \subset \mathring{Q}, \qquad \mathcal{L}^2(A \Delta D) \leq \varepsilon, \qquad \mathcal{P}(D) \leq \mathcal{P}(A) + \varepsilon.$$

PROOF.   For the proof, we refer the reader to [10].   □

In addition to that, in order to create pieces of the interface that are well localized we will need a lower bound for the existence of an interface inside a certain region. Let $x$ be a site in $\mathbb{Z}^2 + (1/2, 1/2)$, and let

$$H(0, nx, \sqrt{n}) = \{y \in \mathbb{R}^2 : d_\infty(y, [(1/2, 1/2), nx]) \leq \sqrt{n}\}$$

be the region formed by the points of $\mathbb{R}^2$ that are at a $\infty$-distance less than $\sqrt{n}$ from the segment $[(1/2, 1/2), nx]$. We will need a lower bound on the following event:

$$\text{Wall}(0, nx, \sqrt{n}) = \{(1/2, 1/2) \leftrightarrow nx \text{ by an open dual path in } H(0, nx, 2\sqrt{n})\}.$$

LEMMA 23.   Let $x \in Q + (1/2, 1/2)$ be a dual site and let $a > 5$. If $p \downarrow p_c$ and $n \uparrow \infty$ in such a way that $n > (p - p_c)^{-a}$, then

$$\liminf_{(n,p)} \frac{1}{(p - p_c)n} \log \Phi^{w,p}_{\Lambda(n)}[\text{Wall}(0, nx, \sqrt{n})] \geq -\tau_c |x|.$$

PROOF.   We consider the case $x = (1,0) + (1/2, 1/2)$, the proof for a general $x$ is similar. Let $M > 0$ be an integer. We denote by $n = Mq + r$ the Euclidean division of $n$ by $M$. Using translation invariance and the FKG inequality, we get that

$$\log \Phi^p_\infty[\text{Wall}(0, nx, \sqrt{n})] \geq q \log \Phi^p_\infty[\text{Wall}(0, Mx, \sqrt{M})]$$
$$+ \log \Phi^p_\infty[\text{Wall}(0, rx, \sqrt{r})].$$

Now, note that the event $\text{Wall}(0, Mx, \sqrt{M})$ is realized as soon as $0 \leftrightarrow Mx$ by an open dual path and there exists no open dual path from $H(0, nx, \sqrt{n}/2)$ to $H(0, nx, \sqrt{n})^c$. Thus, by Proposition 4,

$$\Phi^p_\infty[\text{Wall}(0, Mx, \sqrt{M})] \geq \Phi^p_\infty[0 \leftrightarrow Mx \text{ by an open dual path}]$$
$$- 4n^{3/2} e^{-\lambda(p - p_c)\sqrt{n}},$$

where $\lambda$ is a positive constant. This gives us the following lower bound:

$$(40) \quad \begin{aligned} &\frac{1}{(p - p_c)n} \log \Phi^p_\infty[\text{Wall}(0, nx, \sqrt{n})] \\ &\geq \frac{q}{(p - p_c)n} \log \Phi^p_\infty[0 \leftrightarrow Mx \text{ by an open dual path}] \\ &\quad + \frac{q}{(p - p_c)n} \log\left[1 - \frac{4n^{3/2} \exp(-\lambda(p - p_c)\sqrt{n})}{\Phi^p_\infty[0 \leftrightarrow Mx \text{ by an open dual path}]}\right] \\ &\quad + \frac{1}{(p - p_c)n} \log \Phi^p_\infty[\text{Wall}(0, rx, \sqrt{r})]. \end{aligned}$$

Now we consider a double sequence $M \uparrow \infty, p \downarrow p_c$ such that $M(p - p_c)/\log M \uparrow \infty$, and we suppose also that $(p - p_c)\sqrt{n}/M \uparrow \infty$. Then, by Proposition 4 we have that:

- $\dfrac{q}{(p - p_c)n} \sim \dfrac{1}{(p - p_c)M}$,

- $\lim_{M,p} \dfrac{1}{(p - p_c)M} \log \Phi_\infty^p[0 \leftrightarrow Mx$ by an open dual path$] = -\tau_c$,

- $\lim_{n,M,p} \log \left[ 1 - \dfrac{4n^{3/2} \exp(-\lambda(p - p_c)\sqrt{n})}{\Phi_\infty^p[0 \leftrightarrow Mx \text{ by an open dual path}]} \right] = 0$.

Since $r < M$ we have by the finite energy property that there exist two constants $c, \lambda > 0$ such that $\Phi_\infty^p[\mathrm{Wall}(0, rx, \sqrt{r})] \geq c \exp(-\lambda M)$; thus in the regime specified above we have that

$$\lim_{n,M,p} \frac{1}{(p - p_c)n} \log \Phi_\infty^p[\mathrm{Wall}(0, rx, \sqrt{r})] = 0.$$

Thus, from (40) we get that for every sequence $(n, p) \to (\infty, p_c)$ such that there exists $M$ satisfying $(p - p_c)\sqrt{n}/M \uparrow \infty$ and $M(p - p_c)/\log M \uparrow \infty$, then

$$\liminf_{(n,p)} -\frac{1}{(p - p_c)n} \log \Phi_\infty^p[\mathrm{Wall}(0, nx, \sqrt{n})] \geq \tau_c.$$

The result for the finite volume measure $\Phi_{\Lambda(n)}^{w,p}$ follows from Lemma 8. □

6.2. *Proof of the lower bound.* Now we have all the ingredients to complete the last part of the large deviation principle, namely the proof of the lower bound.

PROPOSITION 24. *Let $a > 5$ and $\nu \in \mathcal{M}(Q)$. If $n \uparrow \infty$ and $\beta \downarrow \beta_c$ in such a way that $n(\beta - \beta_c)^a \uparrow \infty$, then for any weak neighborhood $\mathcal{U}$ of $\nu$*

$$\liminf_{(n,p)} \frac{1}{n(\beta - \beta_c)} \log \mu_{\Lambda(n)}^+[\sigma_n \in \mathcal{U}] \geq -\mathcal{J}(\nu).$$

PROOF. Let $\nu \in \mathcal{M}(Q)$. The statement is not trivial only if $\mathcal{J}(\nu) < +\infty$. In this case, by the definition of the rate function $\mathcal{J}$, there exists a Borel set $A$ of $Q$ such that $\nu$ is the measure with density $-1_A + 1_{Q \setminus A}$ with respect to the Lebesgue measure and $\mathcal{J}(\nu) = \tau_c \mathcal{P}(A)$. Let $\mathcal{U}$ be a weak neighborhood of $\nu$, and let $\varepsilon > 0$. By Lemma 22, there exists a polyhedral set $D$ such that $\overline{D} \subset \overset{\circ}{Q}$, the measure $\widetilde{\nu}$ with density $-1_D + 1_{Q \setminus D}$ with respect to the Lebesgue measure belongs to $\mathcal{U}$ and $\mathcal{P}(D) < \mathcal{P}(A) + \varepsilon$. By definition, the boundary $\partial D$ is a union of $s$ segments included in $\overset{\circ}{Q}$ which we denote by $[a_i, a_{i+1}], 1 \leq i \leq s$. Thus

$$\sum_{1 \leq i \leq s} |a_i - a_{i+1}| \leq \mathcal{P}(A) + \varepsilon.$$

Now we will give a lower bound for the probability that $\sigma_n$ stays in a neighborhood of $\widetilde{\nu}$. Let $f$ be a continuous function on $Q$. For reasons that will be clear a bit later, we choose

$$(41) \qquad \delta = \frac{\theta}{6(1+\theta)(\mathcal{P}(Q) + \mathcal{P}(D))}.$$

In order to evaluate the probability that $|\sigma_n(f) - \widetilde{\nu}(f)| \leq \varepsilon$, we assume that we are in a regime where $n \uparrow \infty$ and $p \downarrow p_c$ in such a way that $n\delta \uparrow \infty$ and rescale the lattice by a factor $L = \lfloor \delta n \rfloor$. Next, we define the sets

$$\underline{A} = \{\underline{y} \in \mathbb{Z}^2 : B_n(\underline{y}) \cap Q \neq \varnothing\}$$

and

$$\underline{E} = \{\underline{y} \in \underline{A} : B_n(\underline{y}) \cap D^c = \varnothing\}, \qquad \underline{F} = \{\underline{y} \in \underline{A} : B_n(\underline{y}) \cap D = \varnothing\}.$$

By choosing $p$ close enough to $p_c$, we can assume that

$$(42) \qquad \mathcal{L}^2\left(Q \setminus \left(\bigcup_{\underline{y} \in \underline{E} \cup \underline{F}} B_n(\underline{y})\right)\right) \leq 6(\mathcal{P}(Q) + \mathcal{P}(D))\delta.$$

Moreover $|\underline{E}| + |\underline{F}| \leq 1/\delta^2$. The set $\underline{E}$ (respectively, $\underline{F}$) can be regarded as a coarse graining of the region inside $D$ (respectively, outside $D$). To evaluate $\sigma_n(f) - \widetilde{\nu}(f)$ we will restrict ourselves to an event $\mathcal{E}$ that gives nice properties to the blocks and such that $\mathcal{E}$ has the right probability of decay. For this let us define the block process, $(Y(\underline{x}), \underline{x} \in \mathbb{Z}^2)$, as the indicator functions of the events

$$R(B'(\underline{x}), L^{1/8}) \cap V(B(\underline{x}), \varepsilon) \cap W(B(\underline{x}), \varepsilon), \qquad \underline{x} \in \mathbb{Z}^2.$$

We choose $\mathcal{E}$ to be the intersection of the events

$$\{Y(\underline{x}) = 1\}, \qquad \underline{x} \in \underline{E} \cup \underline{F}, \qquad \text{Wall}(A_i, A_{i+1}, n, \sqrt{n}), \qquad 1 \leq i \leq s,$$

where $A_i$ is a site on the dual of $\Lambda(n)$ that is closest to $na_i$. Let us evaluate the probability of $\mathcal{E}$

$$(43) \qquad \begin{aligned} \mathbb{P}_n^+(\mathcal{E}) &\geq \mathbb{P}_n^+\left[\bigcap_{\underline{x} \in \underline{E} \cup \underline{F}} \{Y(\underline{x}) = 1\} \,\Big|\, \bigcap_{1 \leq i \leq s} \text{Wall}(A_i, A_{i+1}, n, \sqrt{n})\right] \\ &\quad \times \Phi_{\Lambda(n)}^{w,p}\left[\bigcap_{1 \leq i \leq s} \text{Wall}(A_i, A_{i+1}, n, \sqrt{n})\right]. \end{aligned}$$

First observe that $\bigcap_{1 \leq i \leq s} \text{Wall}(A_i, A_{i+1}, n, \sqrt{n})$ occurs outside the set $\bigcup_{\underline{x} \in \underline{E} \cup \underline{F}} \widetilde{B'}(\underline{x})$ where

$$\widetilde{B'}(\underline{x}) = \{y \in \mathbb{Z}^2 : d(y, \widetilde{B}(\underline{x})) \leq L/10\}.$$

Let $a > 5$. The estimates of Lemmas 9, 10, 11 ensure that, uniformly over the boundary conditions on $B'(\underline{x})$, the probability $\mathbb{P}_n^+[Y(\underline{x}) = 1]$ goes to one when $n \uparrow \infty$ and $p \downarrow p_c$ in such a way that

$$(44) \qquad n(p - p_c)^a \delta(p) \uparrow \infty.$$

Thus by [29] the first factor of (43) goes to 1. On the other hand, by the FKG inequality and by Lemma 23, if $n > (p - p_c)^{-a}$, then

$$\liminf_{n,p} \frac{1}{n(p - p_c)} \log \Phi_{\Lambda(n)}^{w,p} \left[ \bigcap_{1 \le i \le s} \mathrm{Wall}(A_i, A_{i+1}, n, \sqrt{n}) \right]$$

$$\ge - \sum_{1 \le i \le s} |A_{i+1} - A_i| \tau_c.$$

Combining the previous inequalities we get

$$(45) \qquad \liminf_{n,p} \frac{1}{n(p - p_c)} \log \mathbb{P}_n^+[\mathcal{E}] \ge -\tau_c \mathcal{P}(A) - \tau_c \varepsilon.$$

Now we are left with the evaluation of $|\sigma_n(f) - \widetilde{\nu}(f)|$ when $\mathcal{E}$ occurs. Suppose that $\mathcal{E}$ occurs and let $\underline{E}_i, i \in I$ (respectively, $\underline{F}_j, j \in J$) be the connected components of $\underline{E}$ (respectively, $\underline{F}$). For $i \in I$ (respectively, $j \in J$), all the crossing clusters of the good blocks $B(\underline{x}), \underline{x} \in \underline{E}_i$ (respectively, $\underline{x} \in \underline{F}_j$) are connected and belong to one big cluster that we denote by $C_-^i$ (respectively, $C_+^j$). The events $\mathrm{Wall}(A_k, A_{k+1}, n, \sqrt{n}), 1 \le k \le s$ isolate completely the set $\underline{E}$ from $\underline{F}$; thus for every $i \in I$ and $j \in J$, the two clusters $C_-^i$ and $C_+^j$ are disjoint, and, moreover, $C_-^i$ cannot be connected to $\partial \Lambda(n)$. Now suppose for a while that all the clusters $C_-^i, i \in I$, are colored negatively and that all the clusters $C_+^j, j \in J$, are colored positively. We will see later that this restriction does not decrease the probability too much.

Next, we define

$$\mathcal{S} = \bigcup_{\underline{y} \in \underline{E} \cup \underline{F}} \{ C \subset B(\underline{y}) : C \cap \partial^{\mathrm{in}} B(\underline{y}) = \varnothing \}$$

and keeping in mind the suppositions made in the last paragraph, we do the following decomposition:

$$|\sigma_n(f) - \widetilde{\nu}(f)| \le \frac{1}{\theta n^2} \left| \sum_{C \in \mathcal{S}} \sigma(C) \sum_{x \in C} f(x/n) \right|$$

$$+ \frac{1}{\theta n^2} \left| \sum_{\underline{y} \in \underline{E} \cup \underline{F}} \sum_{\substack{x \in B(\underline{y}) \setminus C(\underline{y}) \\ x \leftrightarrow \partial^{\mathrm{in}} B(\underline{y})}} \sigma(x) f(x/n) \right|$$

$$+ \left| \sum_{\underline{y} \in \underline{E} \cup \underline{F}} \left( \frac{\sigma(C(\underline{y}))}{\theta n^2} \sum_{x \in C(\underline{y})} f(x/n) - \int_{B_n(\underline{y})} f(x) \, d\widetilde{\nu}(x) \right) \right|$$

$$+ \left| \int_{Q \setminus \bigcup_{\underline{y} \in \underline{E} \cup \underline{F}} B_n(\underline{y})} f(x)(d\sigma_n(x) - d\widetilde{\nu}(x)) \right|.$$

The second and the third term can be bounded as in the proof of the exponential contiguity by using the properties of the good blocks and the imposed coloring of the clusters $(C_-^i, i \in I)$ and $(C_+^j, j \in J)$. To deal with the fourth term of the last inequality, we use (41) and (42) to get

$$|\sigma_n(f) - \widetilde{\nu}(f)| \leq \frac{1}{\theta n^2} \left| \sum_{C \in \mathcal{S}} \sigma(C) \sum_{x \in C} f\left( \frac{x}{n} \right) \right| + 2\varepsilon \|f\|_\infty + \varepsilon(\varepsilon + 4) \|f\|_\infty.$$

Thus, the estimate burns down to the analysis of the deviations of the first term in the last inequality. Let $\mathcal{E}'$ be the event $\mathcal{E}$ intersected with the color constraint made above. By the same sort of computations as in the proof of the exponential contiguity, we obtain

$$\mathbb{P}_n^+ \left[ \frac{1}{\theta n^2} \left| \sum_{C \in \mathcal{S}} \sigma(C) \sum_{x \in C} f\left( \frac{x}{n} \right) \right| > \varepsilon \|f\|_\infty \Big| \mathcal{E}' \right] \leq \mathbb{P}_n^+ \left[ \frac{1}{|\mathcal{S}|} \left| \sum_{C \in \mathcal{S}} Y_C \right| > \frac{\varepsilon \theta}{(\delta n)^{1/4}} \Big| \mathcal{E}' \right],$$

where $Y_C = \sigma(C) \sum_{x \in C} f(x/n)/(\|f\|_\infty (\delta n)^{1/4})$. Fix $\omega \in \mathcal{E}'$. Observe that under the measure $\mathbb{P}_n^+[\cdot|\omega]$, the random variables $(Y_C, C \in \mathcal{S})$ are independent and take their values in $[-1, 1]$. So we can apply Theorem 7 to control this deviation, and we get

$$\mu_n^{+,\beta}(|\sigma_n(f) - \widetilde{\nu}(f)| \geq \|f\|_\infty (7\varepsilon + \varepsilon^2))$$

$$\geq \frac{1}{2^{|\underline{E}| + |\underline{F}|}} (1 - \exp(-c\varepsilon^2 \theta^{3/2} n^{3/2})) \mathbb{P}_n^+[\mathcal{E}],$$

where $c$ is a positive constant.

We can do the same reasoning with any finite number of continuous functions $f_1, \ldots, f_k$ to get

$$\mu_n^+[\forall l \in \{1, \ldots, k\} \ |\sigma_n(f_l) - \widetilde{\nu}(f_l)| \geq \varepsilon]$$

$$\geq \frac{1}{2^{|\underline{E}| + |\underline{F}|}} (1 - \exp(-c\varepsilon^2 \theta^{3/2} n^{3/2})) \mathbb{P}_n^+[\mathcal{E}].$$

Finally, since $|\underline{E}| + |\underline{F}| \leq 1/\delta^2(p)$ and by (45) we get that if $n \uparrow \infty$ and $p \downarrow p_c$ in such a way that (44) is satisfied, then for every weak neighborhood $\mathcal{U}$ of $\nu$,

$$\forall \varepsilon > 0 \qquad \liminf_{n, \beta} \frac{1}{n(\beta - \beta_c)} \log \mu_n^{+,\beta}[\sigma_n \in \mathcal{U}] \geq -\mathcal{J}(\nu) - \varepsilon \tau_c.$$

Sending $\varepsilon$ to 0 yields the desired lower bound.  $\square$

LABORATOIRE DE MATHÉMATIQUES
UNIVERSITÉ PARIS-SUD
BÂT. 425
91405 ORSAY CEDEX
FRANCE
E-MAIL: rcerf@math.u-psud.fr

64, RUE DE RIVE
1260 NYON
SWITZERLAND
E-MAIL: messikh@gmail.com